\documentclass[%
 aip,
 amsmath,amssymb,
 reprint,%
]{revtex4-1}

\usepackage{graphicx}
\usepackage{dcolumn}
\usepackage{bm}
\pdfoutput=1
\usepackage[dvipsnames]{xcolor}
\usepackage[utf8]{inputenc}
\usepackage[T1]{fontenc}
\usepackage{mathptmx}
\DeclareMathAlphabet{\mathcal}{OMS}{cmsy}{m}{n}
\usepackage{etoolbox}
\usepackage{amsthm}
\usepackage{amsmath}
\usepackage{caption}
\usepackage{subcaption}
\usepackage{physics}
\usepackage{float}

\newtheorem{theorem}{Theorem}[section]
\newtheorem{proposition}[theorem]{Proposition}

\newtheorem{corollary}[theorem]{Corollary}
\newtheorem{definition}[theorem]{Definition}

\newtheorem{remark}[theorem]{Remark}

\newcommand{\T}{\mathfrak{T}}
\newcommand{\tp}{\mathfrak{t}}
\newcommand{\B}{\mathfrak{B}}
\newcommand{\bt}{\mathfrak{b}}
\newcommand{\e}{\mathfrak{e}}

\newcommand{\U}{\mathfrak{U}}
\newcommand{\uf}{\mathfrak{u}}
\newcommand{\C}{\mathfrak{C}}
\newcommand{\cf}{\mathfrak{c}}

\newcommand{\x}{\mathbf{x}}
\newcommand{\y}{\mathbf{y}}
\newcommand{\s}{\mathbf{s}}
\newcommand{\rs}{\mathbf{r}}

\newcommand{\paul}[1]{{\color{black} #1}}

\usepackage{xr}
\makeatletter

\newcommand*{\addFileDependency}[1]{
\typeout{(#1)}
\@addtofilelist{#1}
\IfFileExists{#1}{}{\typeout{No file #1.}}
}\makeatother

\makeatletter
\def\@email#1#2{%
 \endgroup
 \patchcmd{\titleblock@produce}
  {\frontmatter@RRAPformat}
  {\frontmatter@RRAPformat{\produce@RRAP{*#1\href{mailto:#2}{#2}}}\frontmatter@RRAPformat}
  {}{}
}%
\makeatother
\begin{document}

\preprint{AIP/123-QED}


\title[Symmetry breaker governs synchrony patterns in neuronal inspired networks]{Symmetry breaker governs synchrony patterns in neuronal inspired networks}
\author{Anil Kumar}
\email{anikumar@clarkson.edu}
\author{Edmilson Roque dos Santos}%
\altaffiliation[Also at ]{Instituto de Ci\^encias Matem\'aticas e Computa\c{c}\~ao, Universidade de S\~ao Paulo, S\~ao Carlos, Brazil
}
\affiliation{ 
Department of Electrical and Computer Engineering, Clarkson
University, Potsdam, New York, U.S.A.
}%
\affiliation{Clarkson Center for Complex Systems Science, Clarkson
University, Potsdam, New York, U.S.A.
}%

\author{Paul J. Laurienti}
\affiliation{Department of Radiology, Wake Forest University School of Medicine, \paul{Medical Center Blvd}, Winston-Salem, 27101, North Carolina, U.S.A}
\author{Erik Bollt}
\affiliation{ 
Department of Electrical and Computer Engineering, Clarkson
University, Potsdam, New York, U.S.A.
}%
\affiliation{Clarkson Center for Complex Systems Science, Clarkson
University, Potsdam, New York, U.S.A.
}%

\date{\today}

\begin{abstract}
\paul{E}xperiments in the human brain reveal switching between different activity patterns \paul{and functional network organization} over time. Recently, multilayer modeling has been employed \paul{across multiple neurobiological levels} (from spiking networks to brain regions) to unveil novel insights into the emergence and time evolution of synchrony patterns. 
\paul{We consider two layers with the top layer directly coupled to the bottom layer. When isolated, the bottom layer would remain in a specific stable pattern. However, in the presence of the top layer, the network exhibits spatiotemporal switching.} The top layer in combination with the inter-layer coupling
acts as a symmetry breaker, governing the bottom layer and restricting the number of allowed symmetry-induced patterns. This structure allows us to demonstrate the existence and stability of pattern states on the bottom layer, but most remarkably, \paul{it enables} a simple mechanism for switching between patterns based on the unique symmetry-breaking role of the governing layer.
We demonstrate that the symmetry breaker prevents complete synchronization \paul{in the bottom layer, a situation that would not be desirable in a normal functioning brain}. We illustrate our findings using two layers of Hindmarsh-Rose (HR) oscillators, employing the Master Stability function approach in small networks to investigate the switching between patterns.
\end{abstract}

\maketitle

\begin{quotation}
Experimental evidence in neuroscience indicates that the human brain exhibits diverse time-dependent synchrony patterns. Synchrony patterns refer to states in which certain brain regions \paul{exhibit coordinated activity while other regions do not}, and these patterns can swiftly transition in response to external or internal stimuli. Inspired by this rapid switching phenomenon in the brain, we introduce a simple switching mechanism between synchrony patterns in networks. Since the brain is inherently heterogeneous, including different types of nerve cells, chemical compositions, various neuronal pathways, and intercellular interactions, incorporating heterogeneity among network components is crucial. Our work proposes a duplex network, where the bottom layer corresponds to the reference network to undergo transitions between synchrony patterns. We consider symmetry-induced synchrony patterns that are accessible by the bottom layer, so the top layer with inter-layer connections acts as a symmetry breaker governing the emerging symmetry-induced patterns. The validity of the findings is tested through numerical simulations, assessing the linear stability of accessible invariant patterns. Our work highlights the significance of multilayer modeling in neuronal systems. 
\end{quotation}

\section{Introduction} 

Synchrony patterns in networks are found across different areas ranging from physics to neuroscience \cite{Arenas2008,Rodrigues2016}. Rather than exhibiting a fully synchronized cluster, networks often manifest diverse levels of synchronization throughout their structure. For instance, EEG recorded data during unihemispheric sleep in some mammals and birds, in which one half of the brain is more active than the other, resembles distinct domains of synchronized elements \cite{Panaggio_2015,Wang2020}. These synchrony patterns encompass synchronized nodes organized into different forms, including incoherent states, chimera states, cluster states, and complete synchronous states. The specific pattern observed depends on the connectivity structure \cite{Pecora2014,Rodrigues2016,Parastesh2021}. 

\begin{figure*}
\centering
\includegraphics[width = 0.7\textwidth]{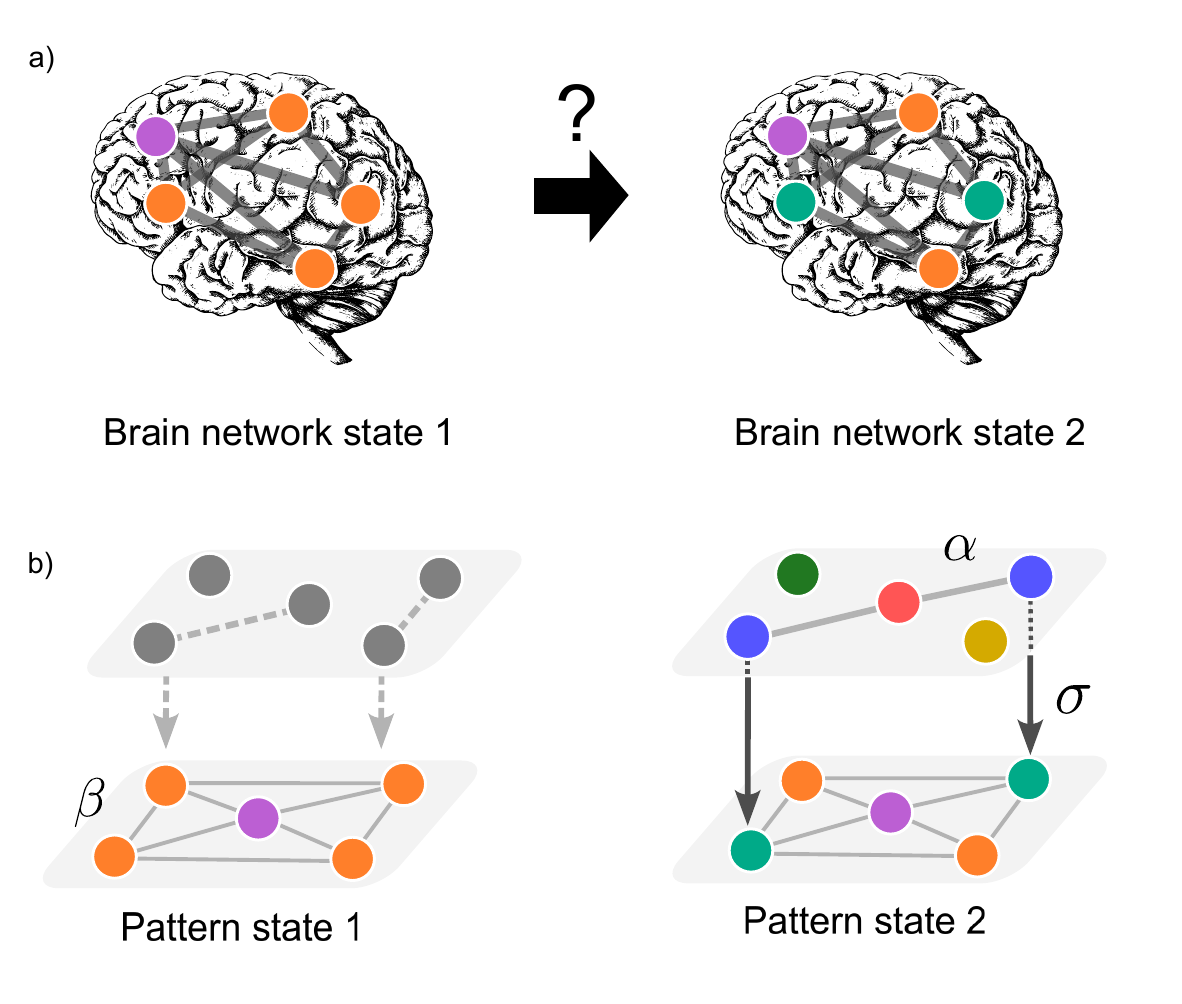}
\caption{\textbf{Symmetry breaker underlies \paul{switching} pattern states in neuronal networks.} a) Schematic diagram of the switching between brain network states. This represents our motivational toy example. b) Each layer contains the same number of nodes, representing different instances of coupled units. Our focus is to observe the change of pattern in the bottom layer by the presence of the top layer and inter-layer connections, the symmetry breaker. The left panel shows the bottom layer exhibits the pattern state 1 for an intra-coupling $\beta$. The right panel shows that once the symmetry breaker is coupled to the bottom layer, the bottom layer shifts its state, exhibiting pattern state 2. Pattern state 1 leaves to be flow-invariant, giving space to pattern state 2. The intra-layer $\alpha$ and inter-coupling $\sigma$ are tuned to pinpoint pattern state 2's stability.}
\label{fig:net_multilayer}
\end{figure*}

Experimental evidence in neuroscience suggests that synchrony patterns are spatiotemporal, evolving dynamically rather than remaining static \cite{buzsaki2006rhythms, Tognoli2014, Deco2015}. Neuroimaging experiments, spanning from functional MRI to EEG measurements of human brain activity, consistently reveal spatiotemporal patterns. \paul{Initial work on functional brain networks were primarily based on an assumption of network stationarity, but a large portion of current research examines dynamic changes in brain networks, even during the resting state \cite{Deco2015, handwerker2012, chang2010, allen2014}.} EEG data recorded during task performance also exhibits a mix of coordination dynamics characterized by both phase locking and metastability \cite{Tognoli2014}. \paul{Thus, consider a simplified} toy example of a brain network depicted in Figure \ref{fig:net_multilayer} a), where the network transitions from one specific \paul{synchronized} pattern to another. This scenario motivates a central question: How do these networks \paul{have} the ability to switch between various synchrony patterns?

In the past few decades, significant effort has been dedicated to understanding synchrony patterns in networks, as evidenced by various studies \cite{Stewart_2003, Golubitsky2005, Belykh_2008}. One well-established criterion to characterize such patterns is network symmetry \cite{GolubitskyMartin2002Tspf, Stewart_2003}. Symmetry refers to a permutation of nodes that leaves the network unchanged, preserving the neighbors of any two nodes that are permuted. The collection of all symmetries in a network forms an algebraic group, which, when applied to coupled identical units, gives rise to clusters. Clusters consist of equivalent nodes under the set of all symmetries, collectively forming a synchrony pattern in the network \cite{Pecora2014}. Symmetry arguments are useful for characterizing synchrony patterns and their stability in networks of identical units \cite{Sorrentino_2016}. However, the brain contains billions of \paul{nerve cells, many different neurotransmitters, various neuronal pathways, and different types and strengths} of interactions between the nerve cells \cite{Koch_1999,bullmore2009}. Therefore, incorporating heterogeneity among network components is crucial for a \paul{more accurate} picture of switching between synchrony patterns in the brain.

Multilayer modeling leverages different features of heterogeneous networks, where distinct node dynamics or coupling functions are considered as layers of a multilayer network \cite{Domenico2013}. In neuroscience, multilayer modeling serves as a valuable tool for representing distinct information derived from the same set of entities and 
{in particular, this kind of model has been important specifically in applications to brain network}  \cite{Crofts_2016,DeDomenico_2017,Battiston_2017,Frolov_2020,Vaiana2020}. 
For instance, Crofts and collaborators \cite{Crofts_2016} investigated important structure-function relations in the Macaque cortical network by modeling it as a duplex network that comprises an anatomical layer, describing the known (macro-scale) network topology of the Macaque monkey, and a functional layer derived from simulated neural activity. Consequently, multilayer modeling has opened the possibility of addressing how neuronal networks switch between synchrony patterns.

Here, we propose a simple mechanism for switching between synchrony patterns in networks using two layers based on the symmetry-breaking role of a governing layer. We consider a network composed of two layers in which layers are named {\it top} and {\it bottom}. The bottom layer serves as the reference network and has several accessible symmetry-induced pattern states. When isolated, the bottom layer maintains a particular pattern state without temporal change. The switching between states emerges due to the presence of the top layer. The top layer and the inter-layer coupling constrains the symmetries of the bottom layer, determining which patterns are permissible --- essentially acting as a \textit{symmetry breaker}. We characterize the existence of symmetry-induced pattern states in a duplex network akin to \cite{DellaRossa2020}. Any bottom layer pattern exists if \paul{it} remains flow-invariant under permutation\paul{,} taking the symmetry breaker into account. The stability of the patterns is assessed in terms of the Master Stability Function approach, which is simplified due to the directionality in inter-layer links. We illustrate the switching between patterns states using numerical simulations in coupled Hindmarsh-Rose oscillators with different node dynamics in each layer as well as different intra\paul{-} and inter-layer coupling functions.

\section{Symmetry breaker: mapping patterns in networks}

We consider that each layer contains $N$ nodes and the evolution of the system is given by:
\begin{align}\label{eq:multilayer_net}
\begin{split}
\dot{x}_i &= \mathfrak{t}(x_i) + \alpha \sum_{j=1}^{N} A_{ij} \uf(x_j), \\
\dot{y}_i &= \mathfrak{b}(y_i) - \beta \sum_{j=1}^{N} L_{ij} \cf(y_j) + \sigma \kappa_{i} D(x_i - y_i), \quad i = 1,2, \dots, N,	
\end{split}
\end{align}
where $x_i, y_i \in \mathbb{R}^n$ are the state vectors of node $i$ in top layer and bottom layer, respectively. We consider the following assumptions:
\begin{itemize}
\item[A.] \textbf{Isolated dynamics.} We consider continuously differentiable top $\tp:\mathbb{R}^n \to \mathbb{R}^n$ and bottom $\bt: \mathbb{R}^n \to \mathbb{R}^n$ isolated dynamics, that 
have an inflowing invariant set, also referred to as dissipative systems. Most smooth nonlinear systems with compact attractors satisfy this assumption. For instance, Hindmarsh-Rose in Equation \eqref{eq:HR} has an attractor \cite{BELYKH_2005,SHILNIKOV_2008}. 

\item[B.] \textbf{Coupling function.} The coupling functions $\uf, \cf: \mathbb{R}^n \to \mathbb{R}^n$ are continuously differentiable. The inter-coupling function is diffusive, depending on the difference of states, and $D \in \mathbb{R}^{n \times n}$.

\item[C.] \textbf{Intra-layer coupling.} The intra-coupling strength is given by $\alpha$ and $\beta$ for the top layer and bottom layer, respectively. The intra-coupling structure of the top layer corresponds to adjacency matrix $A$, where the entries $A_{ij}$ equals $1$ if node $i$ receives a connection from $j$ and $0$ otherwise. The intra-coupling of the bottom layer is given by a Laplacian matrix $L$. For both layers, the connectivity structure is undirected. So, the adjacency matrix and Laplacian matrix are symmetric, and consequently, have real eigenvalues.

\item[D.] \textbf{Directed one-to-one but not onto inter-layer coupling.} The symmetry of the whole multilayer structure moderates what pattern states are possible, and so changes between the layers of the layers allows the top layer to act as a governor over the bottom layer.  We consider a non-symmetric inter-coupling between the two layers (unidirectionally) as illustrated in Figure \ref{fig:net_multilayer}. The top layer drives the bottom layer through a one-to-one but not subjective coupling represented by the diagonal inter-layer matrix 
$K = \mathrm{diag}(\kappa_1, \dots, \kappa_N) \in \mathbb{N}^{N \times N}$ (off-diagonal terms are zero) with each entry $\kappa_i \in \{0, 1\}$, and weighted by an overall inter-coupling strength $\sigma$. For instance, consider a $5 \times 5$ matrix
\begin{align*}
    K = \begin{pmatrix}
        1 & & & & \\
        & 0& & & \\
        & & 1& & \\
        & & & 1& \\
        & & & &0 \\
    \end{pmatrix}.
\end{align*} 
\end{itemize}

Consider an initial state in the network of the five nodes, see the bottom left panel in Figure \ref{fig:net_multilayer}. By the flow-invariance, if the bottom layer starts at pattern state $1$, it remains there and never shifts to pattern state $2$. To allow the bottom layer to attain the other pattern state, we build up the idea of extending the phase space. We assume there exists a top layer, which evolves under different dynamics, coupled through an inter-layer connectivity structure. Out of all possible pattern states permitted in the bottom layer, only a subset is allowed due to symmetry constraints induced by the coupling with the top layer. We view the top layer and inter-layer coupling as a symmetry breaker that restricts, breaks, and changes symmetries on the bottom layer, see Figure \ref{fig:net_multilayer} for an illustration. To attain another pattern state, the symmetry breaker plays the role of mapping one pattern to another. Depending on the structural configuration of the top layer and inter-layer, the symmetry breaker shifts the initial pattern of the bottom layer to another state, which would be impossible in the absence of the symmetry breaker. Figure \ref{fig:net_multilayer} illustrates the bottom layer switching from one pattern to another depending on the symmetry breaker structure.

The following sections are organized as follows. 
\paul{First, in Section \ref{sec:bottom_layer_symm} it is shown that} the topology of the bottom layer allows possible symmetry-induced clusters, and consequently, pattern states. Section \ref{sec:restricting_symmetry} is devoted to characterizing how the flow invariance of Equation \eqref{eq:multilayer_net} and the structure of the symmetry breaker restrict symmetry-induced patterns in the bottom layer. Depending on the symmetry breaker structure, a particular pattern state can cease to exist due to a lack of flow invariance. Once the pattern state is flow-invariant, the bottom layer attains this particular state whenever it is stable. Section \ref{sec:stable_patterns} describes the linear stability of the full system (symmetry breaker and bottom layer). Section \ref{sec:numerical_sim} shows the numerical simulations in coupled Hindmarch-Rose oscillators to demonstrate our findings. 
\paul{Section \ref{sec:dis_con} provides our discussion and conclusions.}

\vspace{0.2cm}
\noindent
\textbf{Notation.} Each vector $u_i \in \mathbb{R}^{n}$ is denoted $u_i = (u_i^1, \dots, u_i^n)$. The vector space $\mathbb{R}^{n}$ is endowed with the $\ell_1$ norm $\|u_i\|_1 = \sum_{k = 1}^n |u_i^k|$. The state space $\mathbb{R}^n \otimes \mathbb{R}^{N}$ can be canonically identified with $(\mathbb{R}^{n})^N$, which we will use for shorter notation. Each element of the space is $\x = \mathbf{vec}(x_1,\dots,x_N) \in (\mathbb{R}^n)^N$, where $\mathbf{vec}$ denotes the vectorization by stacking multiple columns vectors into a single column vector. Also, $(\mathbb{R}^{n})^N$ is equipped with the norm 
\begin{align*}
    \|\x\|_1 = \max_{i \in [N]} \|x_i\|_1.
\end{align*}
Any linear operators on the above spaces will be equipped with the induced operator norm. Finally, let us denote $I_n$ the $n \times n$ identity matrix in $\mathbb{R}^n$.

The notation for the state of the top and bottom layer can be recast as $\x = \mathbf{vec}(x_1,\dots,x_N)$ and $\y = \mathbf{vec}(y_1, \dots, y_N)$, respectively. Moreover, we introduce a tensor representation of Equation \eqref{eq:multilayer_net} into a $2(Nn)$-dimensional system. Define
\begin{align*}
    \T(\x) &= \mathbf{vec}\big(\tp(x_1),\dots, \tp(x_N)\big), \\ 
    \U(\x) &= \mathbf{vec}\big(\uf(x_1),\dots,\uf(x_N)\big), \\
    \B(\y) &= \mathbf{vec}\big(\bt(y_1),\dots, \bt(y_N)\big), \\ 
    \C(\y) &= \mathbf{vec}\big(\cf(y_1),\dots,\cf(y_N)\big).
\end{align*}
Then, Equation \eqref{eq:multilayer_net} can be recast as
\begin{align}
\label{eq:multilayer_tensor_x}
\dot{\x} &= \T(\x) + \alpha (A \otimes I_n) \U(\x) \\ 
\label{eq:multilayer_tensor_y}
\dot{\y} &= \B(\y) - \beta (L \otimes I_n) \C(\y) + \sigma (K \otimes D) (\x - \y),
\end{align}
where $\otimes$ denotes the Kronecker product. 

\section{Bottom layer: Symmetry-induced pattern states}
\label{sec:bottom_layer_symm}
Symmetries of a network are elements of the automorphism group of a graph, acting on the nodes of the network \cite{HararyFrank1969Gt}. Although symmetry is not a necessary condition for synchronization between nodes \cite{Stewart_2003}, it is sufficient for coupled identical oscillators \cite{GolubitskyMartin2002Tspf}. Then, we restrict ourselves to clusters, which are defined as groups of nodes that are synchronized to each other within the same group and distinct from other groups, induced by symmetries of the network, symmetry-based clusters \cite{GolubitskyMartin2002Tspf,Pecora2014}. We will drop the term symmetry-based to characterize the clusters, but it should be clear that we are only dealing with those. More general mechanisms for inducing clusters such as balanced relation \cite{Stewart_2003,Golubitsky2005}, external equitable partition (EEP) \cite{Schaub2016} and graph fibration \cite{DeVille2015} will be explored elsewhere. 

To define cluster, we use orbit partition. The graph automorphism induces a partition of $[N]:= \{1, \dots, N\}$ \cite{HararyFrank1969Gt}. The collection of the graph automorphism orbits acting on $[N]$ induces the orbit partition of $[N]$ \cite{orbit_partition}: $(\bigcup_{l = 1}^{k_{\B}} \mathcal{K}_{\B}^l)$ such that $\mathcal{K}_{\B}^l \cap \mathcal{K}_{\B}^m = \emptyset$ for any $l \neq m \in [k_{\B}]$, where $k_{\B} \in \mathbb{N}$ is the number of clusters and $\mathcal{K}_{\B}^l$ are the clusters. Nodes in $\mathcal{K}_{\B}^l$ can be permuted among each other and will remain synchronized if started in a synchronized state. The trivial orbit partition consists of $\{\{1\},\{2\}, \dots, \{N\}\}$.

From the orbit partition, we define the pattern state. Let $s_l \in \mathbb{R}^n$ be the state of the $l$-th cluster. Then, the \emph{pattern state} of the bottom layer is defined as 
\begin{align}\label{eq:bottom_cluster_state}
 \mathcal{P}_{\B} = \{\y \in \mathcal{V} \subset (\mathbb{R}^n)^N : y_i = s_l \in \mathbb{R}^n, \quad i \in \mathcal{K}_{\B}^l, l \in [k_{\B}]\}
\end{align}
and has dimension $nk_{\B}$, where $\mathcal{V}$ is an inflow invariant set. $\mathcal{P}_{\B}$ corresponds to a manifold embedded in $(\mathbb{R}^{n})^N$. Since any element of the graph automorphism permutes with the Laplacian matrix, the pattern state is an invariant manifold under the bottom layer dynamics, i.e., Equation \eqref{eq:multilayer_net} with $\sigma = 0$. 

\begin{figure*}
\centering
\includegraphics[width = 1.0\textwidth]{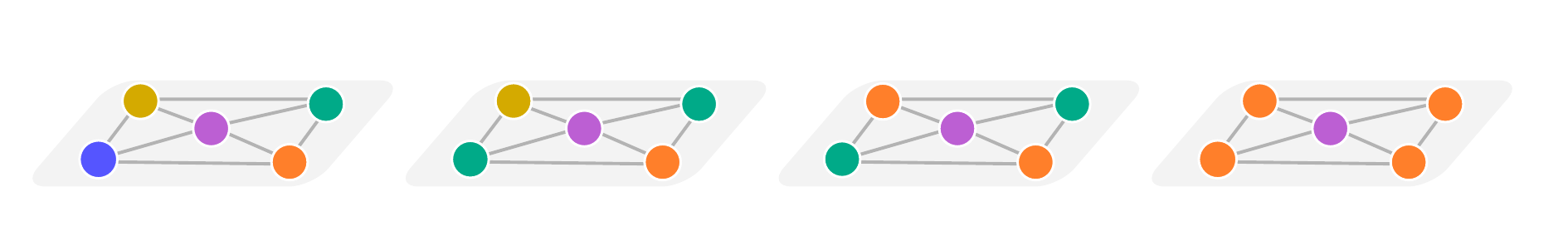}
\caption{\textbf{Symmetry-induced pattern states.} Each orbit partition induced by the graph automorphism and its subgroups creates different patterns in the network with five nodes. Four patterns are illustrated where node color corresponds to a cluster.}
\label{fig:patterns}
\end{figure*}

The graph automorphism has subgroups that also generate different sets of clusters, and consequently, different pattern states \cite{Sorrentino_2016}. To completely characterize all possible pattern states of the bottom layer, we must enumerate the graph automorphism and its subgroups. We employ computational algebra methods/dedicated discrete algebra software \cite{sage} that generate the group automorphism and its subgroups of the bottom layer, as explored in \cite{Pecora2014,Sorrentino_2016,DellaRossa2020}. Figure \ref{fig:patterns} displays different patterns that emerge in a network with five nodes. Each color identifies a different cluster $\mathcal{K}_{\B}^l$ in the network.

\section{Restricting symmetry-induced patterns}
\label{sec:restricting_symmetry}

To define a cluster for a duplex topology, dynamical information must be taken into account: the set of nodes from different layers evolve under different isolated dynamics, see Equation \eqref{eq:multilayer_net}. So, nodes from different layers are not allowed to be permuted between each other. In this case, clusters do not result directly from the orbit partition induced by the graph automorphism of the duplex, but a subgroup where symmetries of the duplex are built up of symmetries that only permute nodes of the same layer \cite{DellaRossa2020}. 

More precisely, let $\mathcal{G}_{\T}$ and $\mathcal{G}_{\B}$ be the graph automorphisms of the top layer and bottom layer, respectively. An element of any of these sets can be represented by a permutation matrix written in a representation in $[N]$. We abuse notation, we only deal with the matrix representatives. Let $P_{\B} \in \mathcal{G}_{\B}$ be a permutation matrix that satisfies $P_{\B} L = L P_{\B}$. In order for Equation \eqref{eq:multilayer_net} to be flow-invariant under the action of this symmetry of the bottom layer, symmetry compatibility \cite{DellaRossa2020} must be achieved, which we state to our specialized scenario: 

\begin{proposition}[Symmetry compatibility \cite{DellaRossa2020}]\label{propo:symm_comp}
Equation \eqref{eq:multilayer_net} is invariant under the action of $P_{\B} \in \mathcal{G}_{\B}$ if there exists $P_{\T} \in \mathcal{G}_{\T}$ such that the conjugacy relation 
\begin{align}\label{eq:conjugacy_relation}
    P_{\B} K = K P_{\T}
\end{align}
is satisfied. 
\end{proposition}
\begin{proof}
To achieve invariance of the Equation \eqref{eq:multilayer_net} under $P_{\B} \in \mathcal{G}_{\B}$, Equations \eqref{eq:multilayer_tensor_x} and \eqref{eq:multilayer_tensor_y} should be invariant under the action of $P_{\B} \otimes I_n$. In fact, denote $\Bar{\x} = (P_{\T} \otimes I_n) \x$ and $\Bar{\y} = (P_{\B} \otimes I_n) \y$. We observe that the dynamics of $\Bar{\x}$ and $\Bar{\y}$ evolve as Equations \eqref{eq:multilayer_tensor_x} and \eqref{eq:multilayer_tensor_y}, i.e., it is invariant, as long as given $P_{\B}$ there exists $P_{\T} \in \mathcal{G}_{\T}$ such that the conjugacy relation is satisfied $P_{\B} K = K P_{\T}$.     
\end{proof}

The conjugacy relation constrains the possible symmetries of the bottom layer, depending on the topology of the top layer and inter-layer structure. The sets
\begin{align}
\begin{split}
    \mathcal{H}_{\T} &= \{P_{\T} \in \mathcal{G}_{\T}: \exists ~ P_{\B} \in \mathcal{G}_{\B} ~ s.t. ~  P_{\B}K = KP_{\T} \}, \\
    \mathcal{H}_{\B} &= \{P_{\B} \in \mathcal{G}_{\B}: \exists ~ P_{\T} \in \mathcal{G}_{\T} ~ s.t. ~ P_{\B}K = KP_{\T} \}
\end{split}
\end{align}
are subgroups of $\mathcal{G}_{\T}$ and $\mathcal{G}_{\B}$, respectively. Since $K$ is not a subjective map, the left and right null spaces are nontrivial. Consequently, we may find more than one $P_{\B}$ that satisfies $P_{\B}K = KP_{\T}$ for a given $P_{\T}$, and vice versa. Then, it is useful to define an equivalence relation between elements in $\mathcal{H}_{\T}$, and similarly, in $\mathcal{H}_{\B}$: 
\begin{itemize}
    \item $P_{\T} \sim_r P_{\T}^{\prime}$ if the right multiplication is equal, $KP_{\T} = KP_{\T}^{\prime}$;
    \item $P_{\B} \sim_l P_{\B}^{\prime}$ if the left multiplication is equal, $P_{\B}K = P_{\B}^{\prime}K$.
\end{itemize}
By the Fundamental Theorem on Equivalence Relations these equivalence relations $\sim_r$ and $\sim_l$ define partitions $\mathcal{H}_{\T}/{\sim_r}$ and $\mathcal{H}_{\B}/{\sim_l}$ on $\mathcal{H}_{\T}$ and $\mathcal{H}_{\B}$, respectively. These partitions are the disjoint union of a finite number of equivalence classes $\{\mathcal{E}_{\T}^i\}_{i}$ and $\{\mathcal{E}_{\B}^i\}_{i}$ as
\begin{align}
    \mathcal{H}_{\T}/{\sim_r} = \bigcup_{i} \mathcal{E}_{\T}^i \quad \mathcal{H}_{\B}/{\sim_l} = \bigcup_{i} \mathcal{E}_{\B}^i.
\end{align}
Then the group of symmetries of the duplex is given by
\begin{align}\label{eq:group_symm_duplex}
    \mathcal{G} = \Big\{\begin{pmatrix}
        P_{\T} & 0 \\ 
        0 & P_{\B} 
    \end{pmatrix} \in \mathbb{R}^{2N \times 2N} : P_{\T} \in \mathcal{E}_{\T}^i, P_{\B} \in \mathcal{E}_{\B}^i\Big\}.
\end{align}
Moreover, we also can define the clusters of the duplex network
\begin{definition}[Duplex clusters]
Let $\mathcal{G}$ be the symmetry group of the duplex. The collection of $\mathcal{G}$-orbits acting diagonally on $[N] \times [N]$ induces the orbit partitions $\{\mathcal{K}_{\T}^j\}_{j = 1}^{k_{\T}}$ and $\{\mathcal{K}_{\B}^l\}_{l = 1}^{k_{\B}}$ with $k_{\T}$ and $k_{\B}$ elements, respectively, that we call clusters of the top and bottom layers, respectively. 
\end{definition}

Note that $N = \sum_{l = 1}^{k_{\T}} |\mathcal{K}_{\T}^l| = \sum_{l = 1}^{k_{\B}} |\mathcal{K}_{\B}^l|$. Also, this definition leads to defining the pattern state of the duplex network. Let $r_j, s_k \in \mathbb{R}^n$ be the state of the $j$-th cluster and $k$-th cluster of the top and bottom layer, respectively. Then, the pattern state of the duplex network is defined as 
\begin{align}
\begin{split}
\mathcal{P} &= \{(\x, \y) \in \mathcal{U} \times \mathcal{V} \subset (\mathbb{R}^n)^N \times (\mathbb{R}^n)^N : x_i = r_l, y_j = s_k \in \mathbb{R}^n,\\
&\qquad i \in \mathcal{K}_{\T}^l, j \in \mathcal{K}_{\B}^k, l \in [k_{\T}], k \in [k_{\B}]\}
\end{split}
\end{align}
with dimension $d:= n(k_{\T} + k_{\B})$, where $\mathcal{U}$ and $\mathcal{V}$ are inflow invariant sets. By construction, $\mathcal{P}$ is an invariant manifold under Equation \eqref{eq:multilayer_net}. Note that from $\mathcal{P}$ definition, it follows that in the absence of the symmetry breaker ($\sigma = 0$), projecting $\mathcal{P}$ onto the second coordinates matches with the pattern state of the bottom layer $\mathcal{P}_{\B}$ in Equation \eqref{eq:bottom_cluster_state}, i.e., $\mathcal{P}_{\B} = \pi_{\B}(\mathcal{P})$, where $\pi_{\B}$ is the canonical projection onto the coordinates of the bottom layer. Once the inter-coupling is positive, this does not necessarily hold, because the particular bottom layer pattern state $\mathcal{P}_{\B}$ may not be invariant under the duplex dynamics in Equation \eqref{eq:multilayer_net}. The symmetry breaker does not only restrict the allowed symmetries of the bottom layer but also constraints on which pattern states are invariant under the duplex dynamics. 

Being invariant under the duplex dynamics also implies other facts that we remark below. 

\noindent
\textbf{Bottom layer clusters: all or nothing.} The inter-layer structure is one-to-one (injective) but not onto (surjective). 
This particular structure implies further constraints on the bottom layer cluster in terms of a dichotomy:
\begin{itemize}
    \item[(i)] Either the bottom layer cluster is not driven by any top layer cluster, 
    \item[(ii)] or it receives one-to-one connections of a top layer cluster that has at least the same size.
\end{itemize}
See Corollary \ref{coro:bottom_layer_dichotomy} for the proof. This result also implies that $k_{\T} \leq k_{\B}$. 

\noindent
\textbf{Symmetry breaker avoids complete synchronization.} The bottom layer complete synchronous state is defined as $y_1 = y_2 = \dots = y_N = s$. In the absence of the symmetry breaker (when $\sigma=0$), this state is an invariant state due to Laplacian coupling. Once the symmetry breaker is present, it loses the flow invariance for any inter-coupling strength $\sigma \neq 0$. In fact, the vector $\mathbf{1} = (1, 1, \dots, 1) \in \mathbb{R}^N$ is an eigenvector of $L$ associated to the eigenvalue $0$, $L \mathbf{1} = 0$. Hence, the intra-coupling term in the bottom layer dynamics Equation \eqref{eq:multilayer_net} cancels out and only the inter-layer coupling remains. The bottom layer's complete synchronous state only remains invariant if two conditions are simultaneously satisfied: the top layer attains the complete synchronous state ($x_1 = x_2 = \dots = x_N$) and the inter-layer coupling matrix is $K = I_N$. Our assumptions violate both conditions, hence the bottom layer's complete synchronous state is not invariant.

\begin{figure*}
 \includegraphics[width=0.7\linewidth]{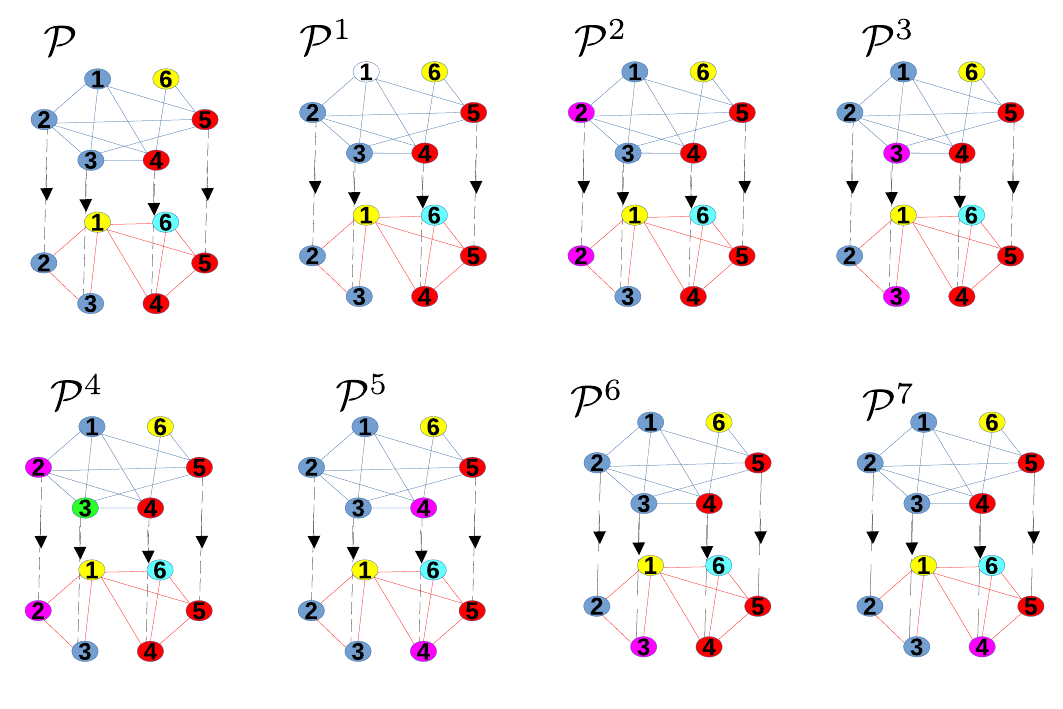}
 \caption{\textbf{Eight symmetry-induced pattern states in a duplex.} Node colors denote the clusters of each layer. The top and bottom layers are illustrated using blue and red edges, respectively. The inter-layer connections are denoted as dashed directed edges. The construction of these patterns can be performed by breaking different synchronous clusters \cite{lodi2021one}. Pattern $\mathcal{P}$ represents the pattern state constructed from the orbital partition of the group in Equation \eqref{eq:group_symm_duplex} while $\mathcal{P}^1-\mathcal{P}^7$ are patterns from breaking different synchronous clusters. $\mathcal{P}^1-\mathcal{P}^4$ are generated from breaking $\mathcal{K}^1_{\T}=\{1,2,3\}$, $\mathcal{P}^5$ from $\mathcal{K}^2_{\T}=\{4,5\}$, $\mathcal{P}^6$ from $\mathcal{K}^2_{\B}=\{2,3\}$, and $\mathcal{P}^7$ from $\mathcal{K}^3_{\B}=\{4,5\}$.}
 \label{fig:irr_construction}
 \end{figure*}

As aforementioned each layer can exhibit multiple admissible symmetry-induced pattern states, see Figure \ref{fig:patterns}. So, from here on we enumerate the pattern states using an additional index. The $i$-th pattern state in the bottom layer is denoted as $\mathcal{P}^i_{\B}$, similarly for the top layer or the duplex. Figure \ref{fig:irr_construction} displays eight distinct symmetry-induced patterns that emerge in the duplex network of each layer containing six nodes. 

\section{Linear stability of pattern states}
\label{sec:stable_patterns}

Flow-invariance under the duplex dynamics is a necessary condition so that the bottom layer may attain any particular pattern state.
This leads to the question of when the pattern state is stable under a certain interval regime of the intra and inter-coupling strengths. In this section, we deduce the linear stability analysis for the duplex dynamics. 

\subsection{Pattern dynamics: quotient dynamics}
\label{sec:quotient_dyn}
The dynamics on the pattern state $\mathcal{P}$ of the duplex network lies in a reduced phase space with dimension $n(k_{\T} + k_{\B})$, where $k_{\T}$ and $k_{\B}$ are the number of clusters in the network. To obtain the equations of motion, also called quotient dynamics, we use the notation in \cite{HahnGena1997GSAM,Schaub2016} to characterize clusters and orbit partition. 

We introduce the notation for the bottom layer, but use the same notation for the top layer, replacing the subscript by $\T$. Let $e_{\B}^l \in \mathbb{R}^N$ be given by the indicator vector that identifies the indices of the nodes in the $l$-th cluster in the bottom layer, i.e., 
\begin{align*}
    (e_{\B}^{l})_i = \begin{cases}
        1, \quad i \in \mathcal{K}_{\B}^l \\ 
        0, \quad \mathrm{otherwise}.
    \end{cases}
\end{align*}
The orbit partition of $[N]$ into $k_{\B}$ clusters is encoded in the characteristic matrix $E_{\B} \in \mathbb{R}^{N \times k_{\B}}$ \cite{HahnGena1997GSAM}: $E_{\B}^{ij} = 1$ if node $i$ belongs to cluster $\mathcal{K}_{\B}^j$ and zero otherwise, i.e.,
the columns of $E_{\B}$ are the indicator vectors $e_{\B}^l$ of the clusters:
\begin{align*}
    E_{\B} = [e_{\B}^1, \dots, e_{\B}^{k_{\B}}]. 
\end{align*}
First note that $E_{\B}^* E_{\B}$ is invertible and diagonal with entries being the size of each cluster $\mathcal{K}_{\B}^l$ \cite{HahnGena1997GSAM}. Also, let $\Pi_{\B} \in \mathbb{R}^{N \times N}$ be defined
\begin{align*}
    \Pi_{\B} = E_{\B} (E_{\B}^* E_{\B})^{-1} E_{\B}^* = E_{\B} E_{\B}^{+},
\end{align*}
where $E_{\B}^{+}$ is the left Moore-Penrose pseudoinverse of $E_{\B}$ and can be seen as a average operator \cite{Schaub2016}. Then, $\Pi_{\B}$ is symmetric, $\Pi_{\B}^2 = \Pi_{\B}$, and $\Pi_{\B}$ and $E_{\B}$ have the same column space. Hence, $\Pi_{\B}$ represents the projection operator onto the column space of $E_{\B}$, i.e., the subspace spanned by the vectors $\{e_{\B}^l\}_{l= 1}^{k_{\B}}$. Moreover, $\Pi_{\B}$ is block-diagonal with each diagonal block a multiple of the `all-ones' matrix \cite{HahnGena1997GSAM} \footnote{The network permutation matrix in \cite{Lin_2016}.}. We denote the projection operator on the column space of $E_{\T}$ and $E_{\B}$ embedded into $\mathbb{R}^{2N}$ as
\begin{align}\label{eq:proj_duplex}
    \Pi = \begin{pmatrix}
        \Pi_{\T} & 0 \\
        0 & \Pi_{\B} 
    \end{pmatrix}.
\end{align}
Denote $\rs = \mathbf{vec}(r_1, \dots, r_{k_{\T}}) \in (\mathbb{R}^{n})^{k_{\T}}$ and $\s = \mathbf{vec}(s_1, \dots, s_{k_{\B}}) \in (\mathbb{R}^{n})^{k_{\B}}$. Then, the isolated dynamics and coupling function evaluated at the pattern state $(\rs, \s)$ of the duplex network are given by 
\begin{align*}
\T_{\rs}(\rs) &= \mathbf{vec}(\tp(r_1), \dots, \tp(r_{k_{\T}})), \\
\U_{\rs}(\rs) &= \mathbf{vec}(\uf(r_1), \dots, \uf(r_{k_{\T}})), \\ 
\B_{\s}(\s) &= \mathbf{vec}(\bt(s_1), \dots, \bt(s_{k_{\B}})), \\
\C_{\s}(\s) &= \mathbf{vec}(\cf(s_1), \dots, \cf(s_{k_{\B}})),
\end{align*}
and satisfy 
\begin{align}\label{eq:quotient_invariance}
\begin{split}
    \T((E_{\T} \otimes I_n)\rs) &= (E_{\T} \otimes I_n) \T_{\rs}(\rs) \\
    \U((E_{\T} \otimes I_n)\rs) &= (E_{\T} \otimes I_n)\U_{\rs}(\rs) \\
    \B((E_{\B} \otimes I_n)\s) &= (E_{\B} \otimes I_n) \B_{\s}(\s) \\
    \C((E_{\B} \otimes I_n)\s) &= (E_{\B} \otimes I_n)\C_{\s}(\s). 
\end{split}
\end{align}
To obtain the quotient dynamics we set $\x = (E_{\T} \otimes I_n)\rs$ and $\y = (E_{\B} \otimes I_n)\s$ and replace in Equations \eqref{eq:multilayer_tensor_x} and \eqref{eq:multilayer_tensor_y}. To recast as a reduced system, we introduce new matrices that are the quotient versions of $A$, $L$, and $K$, i.e., satisfy the following relations:
\begin{align*}
    A E_{\T} &= E_{\T} A_{\rs}, \quad L E_{\B} = E_{\B} L_{\s}, \\ 
    K E_{\T} &= E_{\B} K_{\rs}, \quad K E_{\B} = E_{\B} K_{\s}.  
\end{align*}
The new matrices are given by
\begin{align}\label{eq:quotient_coupling_matrices}
\begin{split}
    A_{\rs} &= E_{\T}^{+} A E_{\T} \in \mathbb{R}^{k_{\T} \times k_{\T}} \\
    L_{\s} &= E_{\B}^{+} L E_{\B} \in \mathbb{R}^{k_{\B} \times k_{\B}} \\
    K_{\rs} &= E_{\B}^+ K E_{\T} \in \mathbb{R}^{k_{\B} \times k_{\T}} \\
    K_{\s} &= E_{\B}^+ K E_{\B} \in \mathbb{R}^{k_{\B} \times k_{\B}},
\end{split}
\end{align}
where the one-to-one inter-layer structure $K = \mathrm{diag}(\kappa_1, \dots, \kappa_N) \in \mathbb{N}^{N \times N}$ implies that $K_{\rs}$ is a $k_{\B} \times k_{\T}$ matrix whose entries are given by
\begin{align*}
    K_{\rs}^{ij} &= \frac{1}{|\mathcal{K}_{\B}^i|} \sum_{l \in \mathcal{K}_{\B}^i \cap \mathcal{K}_{\T}^j} \kappa_l,
\end{align*}
and $K_{\s} = \mathrm{diag}(\frac{1}{|\mathcal{K}_{\B}^1|} \sum_{l \in \mathcal{K}_{\B}^1} \kappa_l, \dots, \frac{1}{|\mathcal{K}_{\B}^{k_{\B}}|} \sum_{l \in \mathcal{K}_{\B}^{k_{\B}}} \kappa_l)$, where $|\cdot|$ corresponds to the number of nodes belonging to the cluster. 

As aforementioned, $\mathcal{P}$ is an invariant manifold under the flow, using \eqref{eq:quotient_invariance} the matrices in Equation \eqref{eq:quotient_coupling_matrices} we obtain the dynamics on $\mathcal{P}$ as given by 
\begin{align}\label{eq:quotient_dyn_cluster}
\begin{split}
    \dot{\rs} &= \T_{\rs}(\rs) + \alpha (A_{\rs} \otimes I_n) \U_{\rs}(\rs), \\
    \dot{\s} &= \B_{\s}(\s) - \beta (L_{\s} \otimes I_n) \C_{\s}(\s) + \sigma (K_{\rs} \otimes D)\rs - \sigma (K_{\s} \otimes D) \s.
\end{split}
\end{align}
Although invariance of the flow in Equation \eqref{eq:multilayer_net} is sufficient to deduce the type of patterns allowed in the bottom layer, knowing if trajectories attain a particular pattern is another issue. It remains to show that there is contraction towards $\mathcal{P}$ transversely. To this end, we obtain equations that govern the dynamics near $\mathcal{P}$. We adapt to our case the exposition in \cite{Pereira_2014}.

\subsection{Linear flow close to the pattern state}
\label{sec:linear_flow}

In order to consider stability of pattern states,
we analyze small perturbations away from $\mathcal{P}$ as
\begin{align}\label{eq:coord_splitting}
    \begin{pmatrix}
        \x \\
        \y
    \end{pmatrix} = \begin{pmatrix}
        E_{\T} \otimes I_n & 0 \\
        0 & E_{\B} \otimes I_n 
    \end{pmatrix}\begin{pmatrix}
        \rs \\
        \s 
    \end{pmatrix} + \xi.
\end{align}
The first term in the sum defines a coordinate on $\mathcal{P}$, and the second term $\xi$ is viewed as a perturbation to the pattern state. The coordinate splitting in Equation \eqref{eq:coord_splitting} is associated with a splitting of $(\mathbb{R}^n)^N \times (\mathbb{R}^n)^N$ as the direct sum of subspaces
\begin{align*}
    (\mathbb{R}^n)^N \times (\mathbb{R}^n)^N = \mathcal{F} \oplus \mathcal{F}^{\perp}
\end{align*}
with associated projections 
\begin{align*}
    \pi_{\mathcal{F}}:(\mathbb{R}^n)^N \times (\mathbb{R}^n)^N \to \mathcal{F} \quad \pi_{\mathcal{F}^{\perp}}:(\mathbb{R}^n)^N \times (\mathbb{R}^n)^N \to \mathcal{F}^{\perp}.
\end{align*}
The subspaces $\mathcal{F}, \mathcal{F}^{\perp} \subset (\mathbb{R}^n)^N \times (\mathbb{R}^n)^N$ are determined by embeddings from $(\mathbb{R}^{n})^{d}$ and $(\mathbb{R}^{n})^{(2N - d)}$, respectively, induced by the duplex group of symmetries $\mathcal{G}$. In fact, $\mathcal{F}$ has a representation in terms of the column space of $E_{\T}$ and $E_{\B}$, respectively. Both matrices encode the cluster information of each layer, and the orthogonal complement of their column space forms a representation of $\mathcal{F}^{\perp}$. More precisely, the following result is valid
\begin{proposition}[Block-diagonalization of coupling matrices.]\label{propo:block_diag_coupling}
Consider $A$ and $L$ the adjacency and Laplacian matrices of Equation \eqref{eq:multilayer_net}, respectively. Let $\Pi_{\T}$ and $\Pi_{\B}$ be orthogonal projection operators onto column spaces associated with clusters in the top and bottom layer, respectively. Consider the set of orthonormal eigenvectors $\{v_{\T}^i\}_{i \in [N]}$ and $\{v_{\B}^i\}_{i \in [N]}$ of $\Pi_{\T}$ and $\Pi_{\B}$, respectively, and denote 
\begin{align*}
    T_{\T} = [v_{\T}^1, \dots, v_{\T}^N] \quad \mathrm{and} \quad T_{\B} = [v_{\B}^1, \dots, v_{\B}^N].
\end{align*}
Then, $B = T_{\T}^{*} A T_{\T}$ and $M = T_{\B}^{*} L T_{\B}$ are block-diagonal and denoted as
\begin{align}\label{eq:block_B}
B = \begin{pmatrix}
    B^{\parallel} & 0\\ 
    0 & B^{\perp}  \\
    \end{pmatrix} \quad \mathrm{and} \quad M = \begin{pmatrix}
    M^{\parallel} & 0\\ 
    0 & M^{\perp}  \\
    \end{pmatrix}
\end{align}
with $B^{\parallel} \in \mathbb{R}^{k_{\T} \times k_{\T}}$, $B^{\perp} \in \mathbb{R}^{(N - k_{\T}) \times (N - k_{\T})}$, $M^{\parallel} \in \mathbb{R}^{k_{\B} \times k_{\B}}$ and $M^{\perp} \in \mathbb{R}^{(N - k_{\B}) \times (N - k_{\B})}$. 
\end{proposition}
\begin{proof}
See Appendix \ref{sec:block_diag_coupling}. 
\end{proof}
First, the above result yields a representation for $\mathcal{F}$ and $\mathcal{F}^{\perp}$:
\begin{align*}
\mathcal{F} &= \mathrm{span} \{v_{\T}^1, \dots, v_{\T}^{k_{\T}}, v_{\B}^1, \dots, v_{\B}^{k_{\B}}\}  \otimes \mathbb{R}^{n}, \\
\mathcal{F}^{\perp} &= \mathrm{span}\{v_{\T}^{k_{\T} + 1}, \dots, v_{\T}^{N}, v_{\B}^{k_{\B} + 1}, \dots, v_{\B}^{N}\} \otimes \mathbb{R}^{n}.
\end{align*}
Second, the range of the adjacency matrix $A$ and the Laplacian matrix $L$ can be split into these two subspaces, i.e., it performs a block-diagonalization and the blocks are denoted with $\parallel$ and $\perp$ symbols. Different approaches are also possible such as irreducible representations (IRR) \cite{Pecora2014}, cluster-based coordinates \cite{Cho_2017}, and simultaneous block diagonalization (SBD) \cite{Zhang_2020}.

Here we consider a pattern state to be stable if the transversal perturbations decay exponentially to zero. To determine the linear stability of a pattern state, we study dynamics close to the pattern state. Let $F_{\B}^l \in \mathbb{R}^{N \times N}$ for each $l \in [k_{\B}]$ be defined as 
\begin{align}\label{eq:F_B}
F_{\B}^l = \mathrm{diag}(e_{\B}^l),    
\end{align}
where $\sum_{l = 1}^{k_{\B}} F_{\B}^l = I_{N}$, and similarly for $F_{\T}^l$. To obtain a linear stability analysis, it suffices the linear flow close to the pattern state given by:
\begin{align*}
\dot{\xi}_{\T}&= \Big( \sum_{l = 1}^{k_{\T}} F_{\T}^{l} \otimes D\T(r_l) + \alpha AF_{\T}^{l} \otimes D \U(r_l) \Big) \xi_{\T}  \\
\dot{\xi}_{\B}&= \Big(\sum_{l = 1}^{k_{\B}} \Big(F_{\B}^{l} \otimes D\B(s_l) - \beta LF_{\B}^{l} \otimes D \C(s_l) \Big) - \sigma (K \otimes D) \Big) \xi_{\B} \\
&\qquad + \sigma (K \otimes D)\xi_{\T},
\end{align*}
written in the coordinates $(\xi_{\T}, \xi_{\B})$ along a curve $(\rs(t), \s(t))$ of the quotient dynamics Equation \eqref{eq:quotient_dyn_cluster}. To obtain a decomposition of the linear flow into the subspaces $\mathcal{F}$ and $\mathcal{F}^{\perp}$, Proposition \ref{propo:block_diag_coupling} introduces a change of coordinates 
\begin{align}\label{eq:coordinate_change}
    \begin{pmatrix}
        \xi_{\T} \\
        \xi_{\B}
    \end{pmatrix} \mapsto \begin{pmatrix}
        T_{\T}^* \otimes I_n & 0 \\
        0 & T_{\B}^* \otimes I_n 
    \end{pmatrix}\begin{pmatrix}
        \xi_{\T} \\
        \xi_{\B}
    \end{pmatrix}, 
\end{align}
which we abuse notation and still denote by $(\xi_{\T}, \xi_{\B})$, that splits the perturbation $(\xi_{\T}, \xi_{\B})$ further into parallel and transversal directions
\begin{align}\label{eq:linear_flow_close}
\begin{split}
\dot{\xi}_{\T}&= \Big( \sum_{l = 1}^{k_{\T}} G_{\T}^{l} \otimes D\T(r_l) + \alpha BG_{\T}^{l} \otimes D \U(r_l) \Big) \xi_{\T} \\
\dot{\xi}_{\B}&=  \Big(\sum_{l = 1}^{k_{\B}} G_{\B}^{l} \otimes D\B(s_l) - \beta MG_{\B}^{l} \otimes D \C(s_l) \Big)\xi_{\B}  \\
&\qquad \sigma (K_{\T} \otimes D) \xi_{\T}  - \sigma (K_{\B} \otimes D) \xi_{\B},    
\end{split}
\end{align}
where $G_{\T}^l = T_{\T}^*F_{\T}^lT_{\T}$, $G_{\B}^l = T_{\B}^*F_{\B}^lT_{\B}$, $K_{\T} =  T_{\B}^*KT_{\T}$ and $K_{\B} =  T_{\B}^*KT_{\B}$. 

The linear stability of any top layer pattern does not depend on the bottom layer dynamics, due to the directed coupling inter-layer structure. Hence, to determine if a duplex pattern is stable, the top layer pattern counterpart must be \paul{stable}. For general intra-coupling functions in both layers and inter-coupling functions, this problem can be assessed numerically, restricting the analysis to the columns associated with the transversal directions of the pattern state, as we detail in Section \ref{sec:transverse_numerically}. 

\begin{remark}[Laplacian coupling case] In the case of Laplacian coupling in the top layer, obtaining a stability condition for the coupling strengths $(\alpha, \beta, \sigma)$ becomes a nontrivial task. For $\sigma = 0$, under specific intra-layer structure \cite{Gambuzza_2019} and linear coupling functions (satisfying special spectral conditions \cite{Pereira_2014}) the patterns are composed of independent clusters, and their corresponding stability can be assessed independently as well. In particular, the transversal directions of the pattern state in Equation \eqref{eq:linear_flow_close} can be analyzed in terms of exponential dichotomy \cite{Pereira_2014} or Milnor stability \cite{Gambuzza_2019}. However, for $\sigma \neq 0$ in Equation \eqref{eq:linear_flow_close}, the parallel and transversal perturbations are coupled to each other, requiring further splitting to separate them. This goes beyond the scope of this paper and will be considered in future work.
\end{remark}   

{\color{black}
\subsection{Master stability approach}
\label{sec:transverse_numerically}

To characterize the stability of all clusters composing a pattern, we view the duplex as a single network composed of nodes, which lie on distinct layers, as having different types.
To find the new coordinate system in Equation \eqref{eq:coordinate_change} numerically, we apply a method for obtaining an irreducible block representation via symmetry breaking of a cluster into smaller clusters \cite{lodi2021one}. This method is suitable for directed networks as in our case when compared to other methods \cite{DellaRossa2020, Brady2021forget}. We numerically solve all transverse perturbations in Equation~\eqref{eq:linear_flow_close}, along with the quotient dynamics given in Equation~\eqref{eq:quotient_dyn_cluster}. The Lyapunov exponents are estimated \cite{Eckmann_1985,Pikovsky_politi_2016}, tracking exponents associated with each cluster composing a pattern state. We denote $\Lambda$ the largest of all Lyapunov exponents of a cluster. The stability of a pattern state for the top or bottom layer is based on the sign values of the largest exponents corresponding to each nontrivial cluster composing the pattern state. See Appendix \ref{appendix:six_nodes} for further details of the stability characterization of pattern states in the case of a duplex composed of six nodes in each layer.

\section{Numerical simulations}
\label{sec:numerical_sim}

In this section, we present numerical simulations illustrating the switching between pattern states through the symmetry breaker. Fixing the topology of the symmetry breaker, we only vary the coupling strengths $(\alpha, \beta, \sigma)$ to induce the switching. The directed inter-layer connections from top to bottom (Equation~\eqref{eq:multilayer_net}) guarantee that the existence of a pattern state in the top layer depends only on $\alpha$, whereas patterns in the bottom layer depend on all coupling strengths. 

\begin{figure*}
 \includegraphics[width=\linewidth]{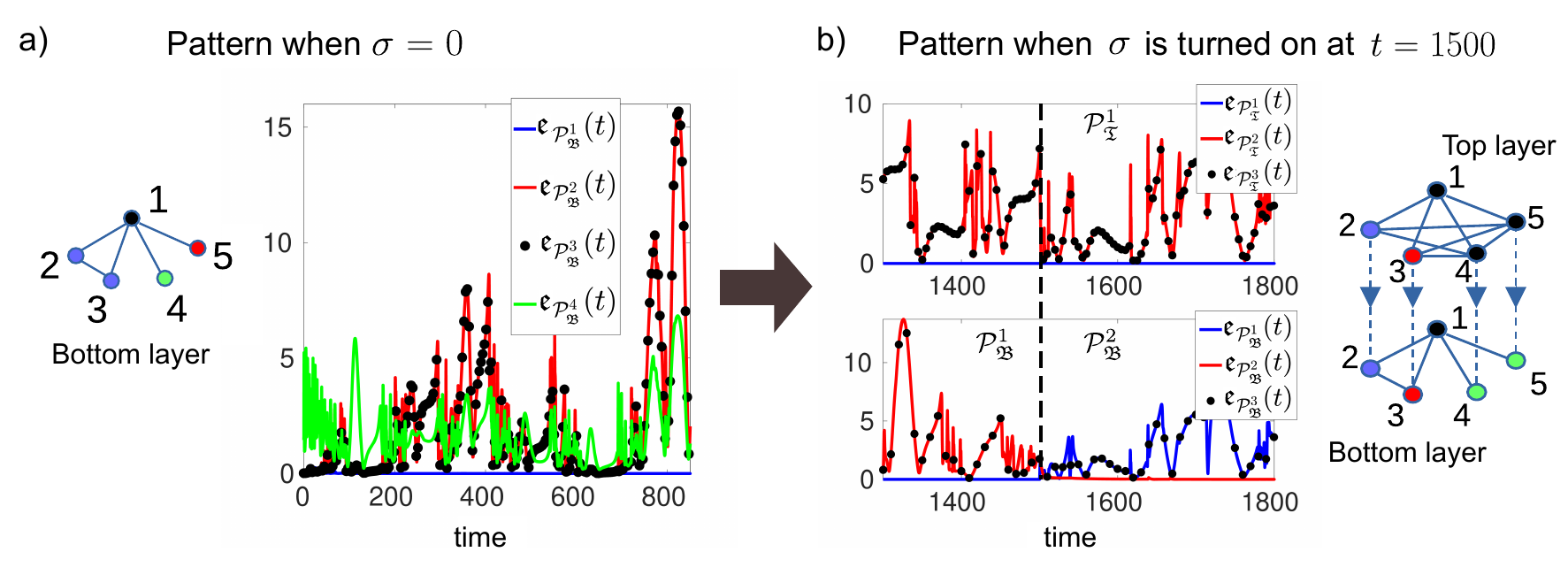}
 \caption{\textbf{Symmetry breaker governs patterns states in coupled HR oscillators.} (a) Initially, when $\sigma=0$, bottom layer is in one pattern state $\mathcal{P}^1_{\B}=\{(a, b, b, c, d)\}$. (b) The bottom layer switches to pattern $\mathcal{P}^2_{\B}=\{(a, b, c, d, d)\}$ in the presence of the symmetry breaker with $\alpha = 0.225$, when $\sigma = 0.5$ at $t=1500$, represented by vertical dashed line. Node color corresponds to nodes in the same cluster, but the same color across layers has no meaning. The HR parameters are $I_\T=3.2, r_\T=0.01, I_\B=3.27, r_\B=0.01, \beta=0.3$. The intra and inter-layer coupling functions are as shown in Equations~\eqref{eq:intra_coupling_fun} and \eqref{eq:inter_coupling_fun}.}
 \label{fig:time_series_plots}
 \end{figure*}

\subsection{Determination of a pattern state} 
\label{sec:synch_error}

We must introduce a measure to quantify when the layer has attained a particular pattern state. Let $\mathcal{K}_{\B}^l$ be the $l$-th cluster in the bottom layer in a pattern $\mathcal{P}_{\B}$. To determine which synchronous pattern exists in the bottom layer, we measure the synchronization error of the pattern $\mathcal{P}_{\B}$, given as  
\begin{align}\label{eq:error_pattern}
    \e_{\mathcal{P}_{\B}}(t) = \frac{1}{k_{\B}} \sum_{l \in [k_{\B}]} \e_{\mathcal{K}_{\B}^l}(t),
\end{align}
where 
\begin{equation}\label{eq:error_cluster}
\e_{\mathcal{K}_{\B}^l}(t) = \frac{1}{|\mathcal{K}_{\B}^l|(|\mathcal{K}_{\B}^l| - 1)} \sum_{i, j \in \mathcal{K}_{\B}^l} \|y_i(t) - y_j(t)\|_1.
\end{equation}
When $\e_{\mathcal{P}_{\B}}(t)$ is \paul{away from zero the pattern does not exist at time $t$,} while $\e_{\mathcal{P}_{\B}}(t)$ close to zero shows that the layer attained $\mathcal{P}_{\B}$. Note that Equation~\eqref{eq:error_cluster} is only defined for nontrivial clusters since trivial clusters have $|\mathcal{K}_{\B}^l| = 1$. Therefore a pattern state is decided by calculating the synchronization errors for all nontrivial clusters in it. Similarly, we can define the synchronization error for the pattern $\mathcal{P}_{\T}$ in the top layer.  

The definition of the synchronization error in Equation~\eqref{eq:error_pattern} allows us to deduce the following observations:
\begin{itemize}
\item [(i)] In case the synchronization error is away from zero for all patterns that contain at least one cluster, then the bottom layer is in the incoherent state, which we denote as $\mathcal{P}^0_{\B}$. 
\item [(ii)] Consider two distinct patterns $\mathcal{P}^i_{\B}$ and $\mathcal{P}^j_{\B}$. Note that whenever $\e_{\mathcal{P}^i_{\B}}(t)$ is close to zero but $\e_{\mathcal{P}^j_{\B}}(t) \ne 0$ as $t \to \infty$, we say the bottom layer is in pattern state $\mathcal{P}^i_{\B}$.
\item[(iii)] If a pattern state $\mathcal{P}^i_{\B}$ contains more nontrivial clusters than $\mathcal{P}^j_{\B}$ and if $e_{\mathcal{P}^i_{\B}}$ and $e_{\mathcal{P}^j_{\B}}$ are both close to zero, the bottom layer is in pattern $\mathcal{P}^i_{\B}$. 
\end{itemize}

Note that the same information about a pattern state can be deduced by plotting $\e_{\mathcal{K}_{\B}^l}(t)$ for all nontrivial clusters $\mathcal{K}_{\B}^l$ as well. For networks exhibiting a small number of pattern states, plotting all synchronization errors corresponding to each nontrivial pattern is a straightforward way to determine which pattern state the network has attained.

\subsection{Hindmarsh-Rose oscillators and coupling functions} 
\label{sec:model}
We focus our attention on a neurologically relevant dynamical system and consider Hindmarsh-Rose (HR) oscillators \cite{Hidmarsh_Rose1984} as nodes in both layers. The bottom isolated dynamics $\bt$ is given in arbitrary coordinates $(v, w, z) \in \mathbb{R}^3$ as
\begin{equation}\label{eq:HR}
(v, w, z) \mapsto 
\begin{bmatrix}
    w -a_{\B} v^3 + b_{\B} v^2 - z +I_{\B}\\
     c_{\B} - d_{\B} v^2 - w\\
    r_{\B} \big(s_{\B} [v - t_{\B}] - z\big)
\end{bmatrix},
\end{equation}
where $a_{\B}, b_{\B}, c_{\B}, d_{\B}, I_{\B}, r_{\B}, s_{\B}, t_{\B}$ are parameters that determine the dynamical regime such as a fixed point, periodic orbit or chaotic dynamics. We fix some of the parameters in Equation~\eqref{eq:HR} as $a_{\B} = 1, b_{\B} = 3, c_{\B} = 1, d_{\B} = 5, s_{\B} = 4, t_{\B} = -0.5(1+\sqrt{5})$\cite{Hidmarsh_Rose1984}. The isolated dynamics in the top layer $\tp$ is the same as $\bt$, except that one or both parameters $I_{\T} \ne I_{\B}$ and $r_{\T} \ne r_{\B}$ are made different so that nodes in the layers are non-identical. 

The intralayer coupling functions are as follows 
\begin{equation}\label{eq:intra_coupling_fun}
   	\uf(x_j)= 
	\begin{bmatrix}
		x^1_j \\
		0\\
		0
	\end{bmatrix}, \;
       \cf(y_j)= 
	\begin{bmatrix}
		y^1_j \\
		0\\
		0
	\end{bmatrix},
\end{equation}
while the inter-layer coupling function follows the relation
\begin{equation}\label{eq:inter_coupling_fun}
	D= 
	\begin{bmatrix}
		0 & 0 & 0 \\
		0 & 1 & 0 \\
		0 & 0 & 0
	\end{bmatrix}.
\end{equation}
Equation \ref{eq:multilayer_net} is solved numerically using the Runge-Kutta fourth-order method with varying integration time steps discarding a transient time so the oscillators reach a steady state. To analyze any clusters $\mathcal{K}_{\T}^l$ and $\mathcal{K}_{\B}^m$ in the top and bottom layers, respectively, the initial conditions are selected such that $x_i(0) \approx x_j(0) ~\forall i,j \in \mathcal{K}_{\T}^l$ and $y_i(0) \approx y_j(0)~ \forall i,j \in \mathcal{K}_{\B}^m$.

\subsection{Switching patterns states}


Our numerical simulations show that the symmetry breaker can drive synchrony patterns in the bottom layer, using $\alpha$ and $\sigma$ as control parameters. The switching between two patterns in the bottom layer can be described as follows: the desynchronization of a cluster in the top layer also breaks the corresponding compatible cluster in the bottom. Likewise, if a synchronized cluster exists in the top layer, increasing diffusive coupling through the $y$ variable eventually synchronizes a compatible cluster in the bottom \cite{huang2009generic}. Therefore, synchronization and desynchronization of clusters in the top layer change patterns at the bottom layer. 

Figure~\ref{fig:time_series_plots} displays the switching between pattern states in the bottom layer of $5$ nodes. When the symmetry breaker is not present, $\sigma = 0$, the bottom layer is in pattern $\mathcal{P}^1_{\B}=\{(a, b, b, c, d)\}$ as shown by the time-evolution of the synchronization error $\e_{\mathcal{P}^1_{\B}}(t)$ close to zero in Figure~\ref{fig:time_series_plots} a) in comparison to the others admissible patterns, $\mathcal{P}^2_{\B}=\{(a, b, c, d, d)\}$ and $\mathcal{P}^3_{\B}=\{(a, b, b, c, c)\}$. Then, for $\sigma \neq 0$ at $t=1500$, we observe that the bottom layer switches from $\mathcal{P}^1_{\B}$ to $\mathcal{P}^2_{\B}$, as illustrated in Figure~\ref{fig:time_series_plots} b). The synchronization error $\e_{\mathcal{P}^1_{\B}}(t)$ starts to oscillate away from zero, whereas $\e_{\mathcal{P}^2_{\B}}(t)$ decays to zero. To achieve this particular switching, the top layer attains the pattern state $\mathcal{P}^1_{\T}=\{(a, b, c, a, a)\}$ as quantified by the synchronization error $\e_{\mathcal{P}^1_{\T}}(t)$, which is constant and close to zero, illustrated in Figure~\ref{fig:time_series_plots} b). The initial conditions of each layer were selected such that the top and bottom layers were close to patterns $\mathcal{P}^3_{\mathfrak{T}}$ and $\mathcal{P}^3_{\mathfrak{B}}$, respectively. Choosing the same HR parameters and initial conditions from Figure~\ref{fig:time_series_plots}, but tuning $\alpha$ parameter, we can also observe that the bottom layer switches from the same pattern $\mathcal{P}^1_{\B}$ to other patterns, see Figure~\ref{fig:other_patterns}. In Figure~\ref{fig:other_patterns} a), the top layer is in the incoherent pattern, inducing that the bottom layer switches to the incoherent pattern $\mathcal{P}^0_{\B} = \{(a, b, c, d, e)\}$, as observed by all synchronization errors being away from zero after $t = 1500$. In Figure~\ref{fig:other_patterns} b), the bottom layer goes to $\mathcal{P}^3_{\B}$, which contains the minimum number of symmetry-induced clusters, as illustrated by the synchronization errors of all patterns decaying to zero once $\sigma$ is turned on. 

\begin{figure}
\includegraphics[width=\linewidth]{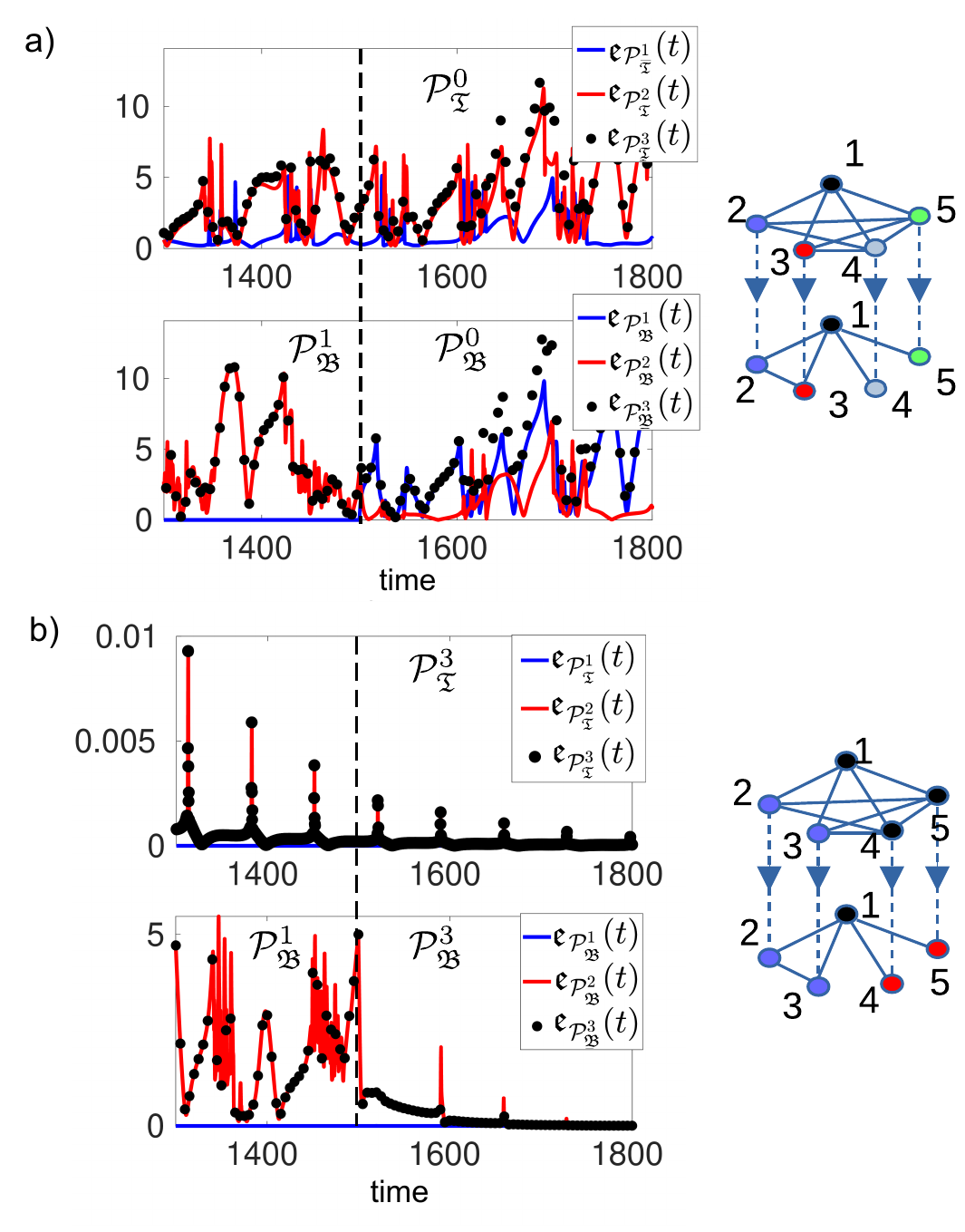}
\caption{\textbf{Bottom layer can attain other patterns by changing $\alpha$.} (a) The bottom layer switches from $\mathcal{P}^1_{\B}=\{(a, b, b, c, d)\}$ to an incoherent pattern $\mathcal{P}^0_{\B} = \{(a, b, c, d, e)\}$ when $\alpha = 0.1$ and $\sigma$ is turned on at $t=1500$. (b) The bottom layer switches from $\mathcal{P}^1_{\B}=\{(a, b, b, c, d)\}$ to $\mathcal{P}^3_{\B}=\{(a, b, b, c, c)\}$ when $\alpha = 0.425$.}
\label{fig:other_patterns}
 \end{figure}

All the switching phenomena occurred from the same pattern state, $\mathcal{P}^1_{\B}$, and the bottom layer could attain any other pattern state apart from the complete synchronous pattern. In other words, from pattern $\mathcal{P}^1_{\B}$ we found three possible transition pathways to other pattern states. But a natural question is: what are the admissible transition pathways starting from a different pattern? Figure~\ref{fig:transition_pathways} displays all admissible transition pathways that switch the bottom layer's pattern. Fixing the topology of the top layer and inter-coupling, the switching occurs only by changing $(\alpha, \sigma)$. We observe that all transition pathways among pattern states form a complete directed graph, where the nodes are the pattern states, and the directed edges are the transition pathways. This confirms that the symmetry breaker can drive the bottom layer to a different pattern state, regardless of which state the bottom layer starts at. 

\begin{figure*}
 \includegraphics[width=0.8\linewidth]{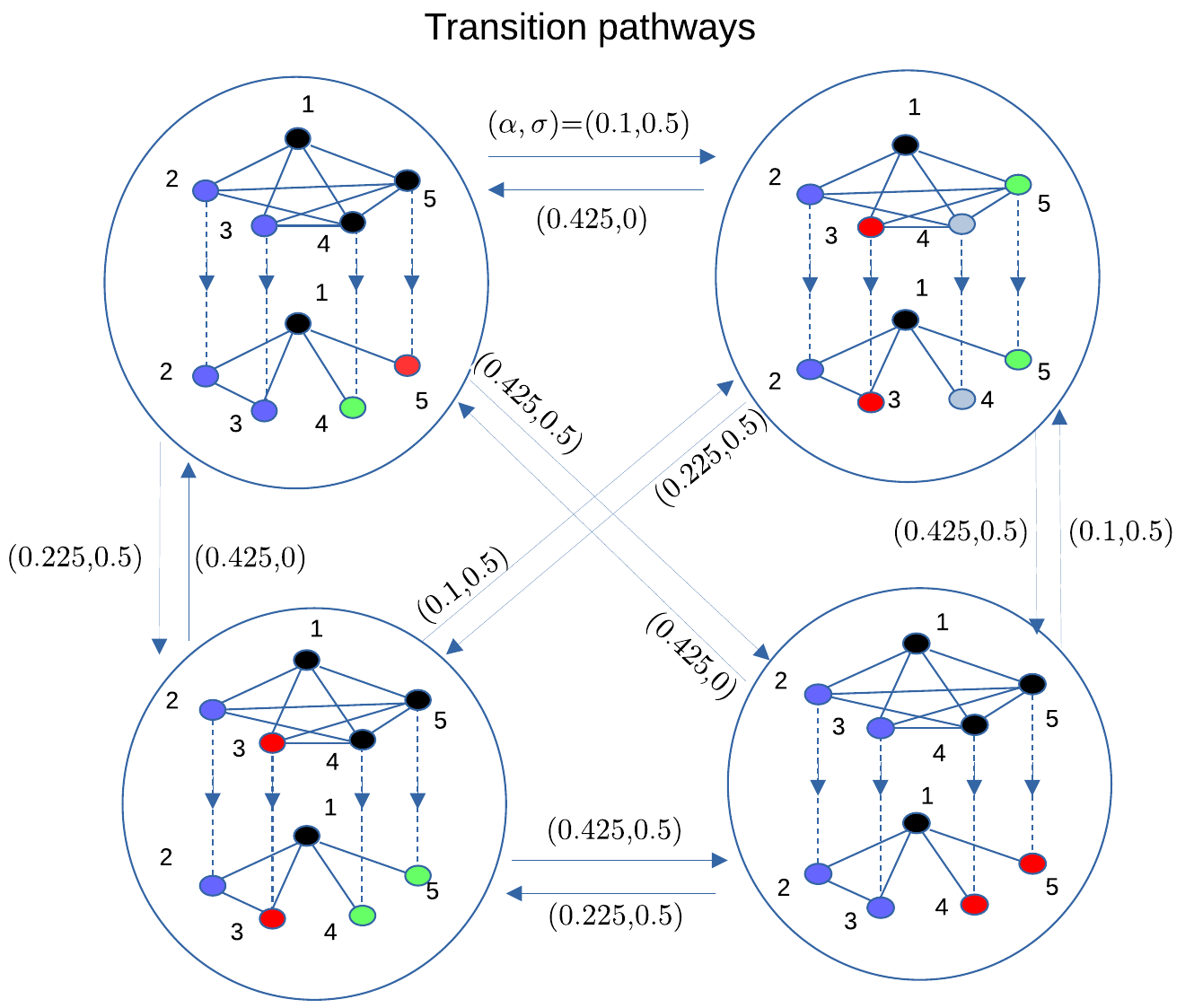}
 \caption{\textbf{Transition pathways driven by symmetry breaker}. For a fixed symmetry breaker topology, regardless of the initial pattern state, the symmetry breaker drives the bottom layer to any other pattern state. There are total of four pattern states: $\mathcal{P}^1_{\B}=\{(a, b, b, c, d)\}$ (left top corner), $\mathcal{P}^0_{\B} = \{(a, b, c, d, e)\}$ (right top corner), $\mathcal{P}^3_{\B}=\{(a, b, b, c, c)\}$ (right bottom corner), and $\mathcal{P}^2_{\B}=\{(a, b, c, d, d)\}$ (left bottom corner). The numerical simulations to draw these transition pathways have the same network dynamics parameters as in Figure~\ref{fig:time_series_plots}.}
 \label{fig:transition_pathways}
\end{figure*}

\section{Discussion and conclusion}
\label{sec:dis_con}

Building upon a multilayer perspective for neuronal network dynamics \cite{Majhi2016, Majhi_2017,omelchenko2019control,mikhaylenko2019weak,zakharova2020chimera,Parastesh2021}, the current work has proposed a mechanism for a network to switch between pattern states. The network is composed of two layers, a duplex network, where the bottom layer is the reference network. The top layer together with the inter-coupling connections forms a symmetry breaker, driving the invariant and stable symmetry-induced patterns in the bottom layer. Instead of characterizing the existence and stability of individual clusters \cite{blaha2019cluster, DellaRossa2020,lodi2021one,sorrentino2020group}, the relevant information is the collection of clusters that defines each symmetry-induced pattern. We have characterized all symmetry-induced pattern states in the bottom layer that are compatible in terms of the symmetry constraints imposed by the connections with the top layer. In particular, our work has demonstrated that the symmetry breaker avoids the complete synchronous state, which corresponds to excessive (abnormal) synchrony and is associated with neurological disorders \cite{Jiruska2013epilepsy}. \paul{Although whole-brain synchrony is not desirable, clusters of synchronized regions that dynamically change membership and activity patterns, as shown here, likely underlie complex cognitive and emotional processes \cite{buzsaki2006rhythms, pessoa2014understanding}.}

Our work has mapped the transition pathways of a set of admissible pattern states for small networks, fixing the symmetry breaker's topology. However, symmetry breakers are also capable of inducing specific and targeted symmetry-induced pattern states in networks. Hence, instead of fixing the topology, an interesting research direction is a control perspective \cite{omelchenko2019control,Gabuzza_2021}: from a set of chosen symmetry-induced patterns in the bottom layer, to design the topology of the symmetry breaker such that the bottom layer can attain these particular pattern states. Although the inter-coupling connections are directed, which simplifies the analysis, similar results are also valid for bidirectional inter-coupling connections.

Our results will be explored in other network dynamics that include neurologically relevant information such as brain-inspired network topology and bio-physical details in the network dynamics. \paul{When building more neurobiologically relevant models, it will be essential to consider how the multilayer network may manifest in the brain. The multilayer organization can be literally true with a top layer controlling the bottom layer. This structure is unlikely in a distributed system like the brain. It is also possible that the multilayer structure is analogical with all nodes in the same layer but having unique roles. Such an organization is much more likely in the brain with some nodes serving the symmetry-breaker role and controlling the pattern states in other nodes. One could even imagine that this symmetry-breaker role is dynamic with nodes flipping from the ``top layer'' to the ``bottom layer'' and vice versa depending on the state of the brain and the incoming sensory stimuli.} For small networks, \paul{like the models used here}, the pattern states of the bottom layer can be accessed numerically via the time-evolution of the synchronization error \eqref{eq:error_pattern}. \paul{However, characterizing all admissible for larger networks, like neurobiologically-relevant ones, becomes a nontrivial task}. So, numerical techniques to identify different pattern states via clustering \cite{bollt2023fractal} seems a promising direction. In summary, our work illustrates a mechanism by which one network can assist or drive patterns of synchrony in others, \paul{providing} potential insights into switching between patterns in a complex system such as the brain.

\begin{acknowledgments}
 A.K., E.R.S., P.J.L., and E.B. acknowledge support from \paul{the Collaborative Research in Computational Neuroscience (CRCNS) through R01-AA029926.} E.B. was also supported by ONR, ARO, DARPA RSDN, and AFSOR. E.R.S. was also supported by Serrapilheira Institute (Grant No. Serra-1709-16124). The brain graphic in Figure \ref{fig:net_multilayer} was extracted from FreePick.
\end{acknowledgments}

\section{Author Declarations}

\noindent
\textbf{Conflict of Interest.}
The authors have no conflicts to disclose.

\noindent
\textbf{Author Contributions.}
\textbf{Anil Kumar}: Conceptualization (equal); Investigation (equal); Software (lead); Validation (equal); Visualization (equal); Writing – original draft (equal); Writing – review \& editing (equal). \textbf{Edmilson Roque dos Santos}:  Conceptualization (equal); Formal analysis (lead); Investigation (equal); Validation (equal); Visualization (equal); Writing – original draft (equal); Writing – review \& editing (equal). \paul{\textbf{Paul Laurienti}: Funding acquisition (equal); Writing – review \& editing (equal).} \textbf{Erik Bollt}: Conceptualization (equal); Funding acquisition (equal); Investigation (equal); Methodology (lead); Supervision (lead); Writing – review \& editing (equal).

\section*{Data Availability}
The data that support the findings of this study are available from the corresponding author upon reasonable request.

\appendix

\section{Appendixes}

\subsection{Bottom layer's cluster: all or nothing}
The special choice of inter-coupling structure implies the following result:
\begin{corollary}[Bottom layer's cluster: all or nothing.]\label{coro:bottom_layer_dichotomy}
Consider a bottom layer cluster of size $m$. Then, one of the two scenarios holds: 
\begin{itemize}
\item[(i)] The bottom layer cluster does not receive connections at all. 
\item[(ii)] The bottom layer cluster receives one-to-one connections from a top layer cluster, which has at least the same size $m$.
\end{itemize}
\end{corollary}
The claim follows from the flow-invariance of Equation \eqref{eq:multilayer_net}, and consequently, from the symmetry compatibility \ref{propo:symm_comp}. 
\begin{proof}
Consider top and bottom layers clusters $\mathcal{K}_{\T}$ and $\mathcal{K}_{\B}$, respectively. Denote $m$ as the size of the bottom layer cluster, i.e., $m = |\mathcal{K}_{\B}|$. Let $K = \mathrm{diag}(\kappa_1, \dots, \kappa_N)$ be the inter-layer coupling.

Note that if node $i$ and $j$ are in cluster $\mathcal{K}_{\B}$, then there exists a permutation matrix $P_{\B}$ that permutes these nodes and satisfies the symmetry compatibility in Equation \eqref{eq:conjugacy_relation}. So, there exists a permutation matrix $P_{\T}$ that 
\begin{align}\label{eq:P_star}
 K P_{\T} = P_{\B} K.   
\end{align}
Also, note that $P_{\B}$ permutes rows $i$ and $j$ of $K$. Since $K$ is a diagonal matrix, to attain the Equation \eqref{eq:P_star}, $P_{\T}$ must permute columns $i$ and $j$ of $K$. This implies that the entries of $K$ satisfy $\kappa_i = \kappa_j$ for $i, j \in \mathcal{K}_{\T}$. Repeating the argument to every node in the bottom layer $\mathcal{K}_{\B}$, we obtain a constraint over the entries of $K$:
\begin{align*}
    \kappa_{i_1} = \kappa_{i_2} = \dots = \kappa_{i_m} = \begin{cases}
        1, \\
        0
    \end{cases}, \quad i_l \in \mathcal{K}_{\B}.
\end{align*}
Evaluating both cases implies the claim.     
\end{proof}

\begin{figure*}[t]
\centering
  \includegraphics[width=0.8\textwidth]{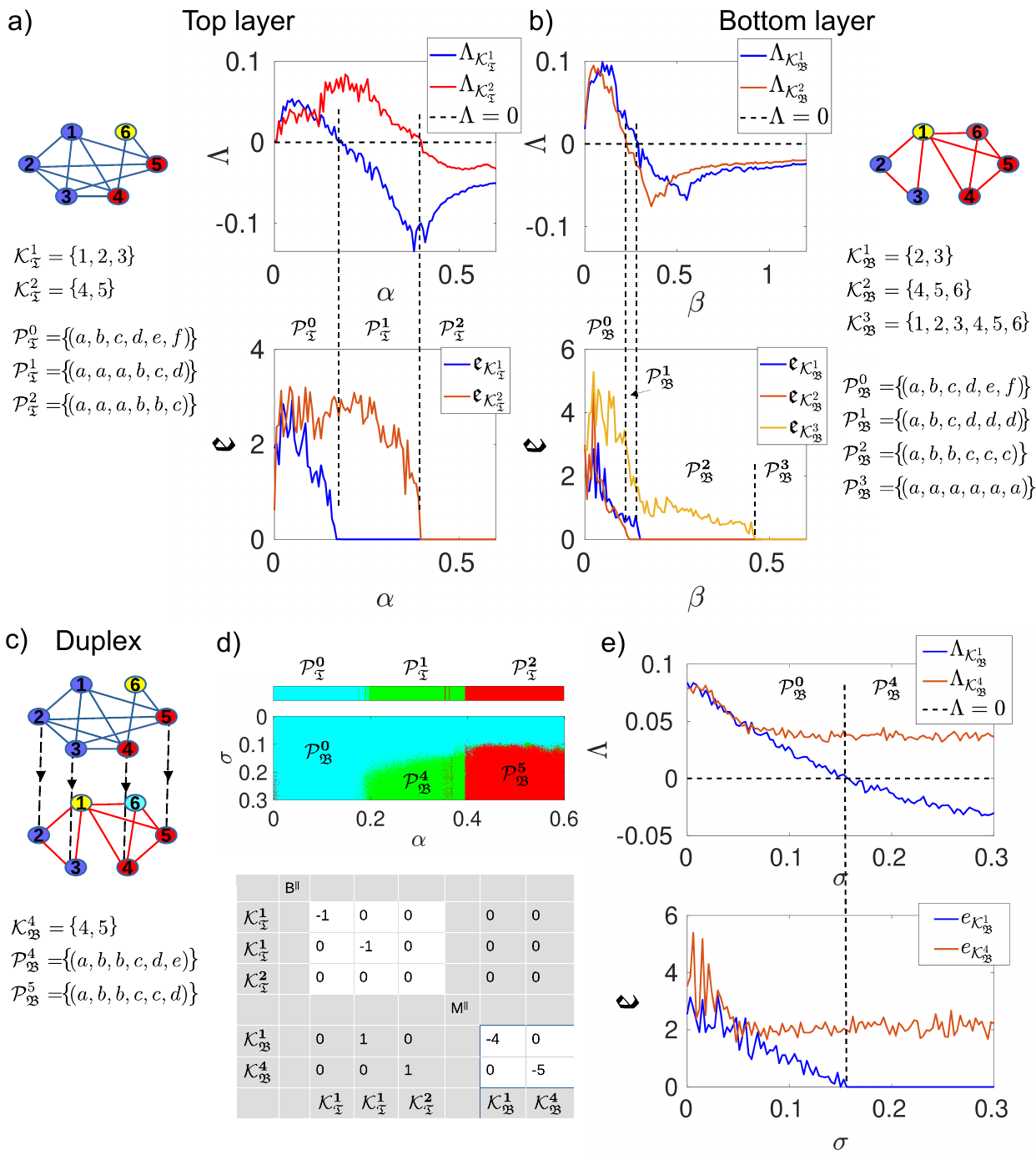}
\caption{\textbf{Stability characterization for switching between patterns} (a) Pattern states in the top layer are depicted by the change in the largest Lyapunov exponent corresponding to each nontrivial cluster, and a comparison with synchronization errors from numerical simulations. (b) Pattern states in the bottom layer when $\sigma=0$. (c) A multilayer network is constructed by joining the top and bottom layers through directed inter-layer links. (d) The top panel shows the stability analysis of the pattern states in the duplex network in $(\alpha, \sigma)$ space. The colors represent when the associated largest Lyapunov exponent is negative for the particular pattern. The bottom panel shows the transformed coupling matrix in the node space. (e) The stability analysis and time-averaged synchronization errors for compatible clusters in the bottom layer as $\sigma$ changes. Here we choose $\alpha = 0.2950$ (e) and $\beta=0.05$ ((d), (e)). The HR parameters in all figures are are $I_{\T}=3.1, r_{\T}=0.01, I_{\B}=3.2, r_{\B}=0.005$ with coupling functions given by Equations~\eqref{eq:intra_coupling_fun}, \eqref{eq:inter_coupling_fun}.}
\label{fig:directed_multilayer_n_12}
\end{figure*}

\subsection{Proof of Proposition \ref{propo:block_diag_coupling}}
\label{sec:block_diag_coupling}
\begin{proof}
We develop the proof for the Laplacian matrix, but the arguments can be repeated to the adjacency matrix. We split the proof into two steps.

\noindent
\textbf{Eigenspace of $\Pi_{\B}$.} Since $\Pi_{\B}$ is an orthogonal projection operator onto the column space of $E_{\B}$. Then, $\{e_{\B}^l\}_{l = 1}^{k_{\B}}$ is the set of eigenvectors associated to the eigenvalue $1$. So, denote $v_{\B}^l = e_{\B}^l$ for $l = 1, \dots, k_{\B}$. 

Consider the orthogonal complement of $\mathrm{span}\{e_{\B}^l\}_{l = 1}^{k_{\B}}$, where any vector is in the kernel of $\Pi_{\B}$. To obtain a set of orthonormal eigenvectors $\{v_{\B}^i\}_{i \in [N]}$ which diagonalizes $\Pi_{\B}$, it suffices to extend the basis $\{v_{\B}^l\}_l$ to a basis of $\mathbb{R}^{N}$, and apply the Gram Schmidt process.

\noindent
\textbf{$L$ and $\Pi_{\B}$ commute.} First, note that $\Pi_{\B}$ is a doubly stochastic matrix. In fact, let $\mathbf{1}$ be the all ones vector, then $\Pi_{\B} \mathbf{1} = \mathbf{1}$, and $\Pi_{\B}$ is symmetric, so $\Pi_{\B}$ is doubly stochastic. Using Birkhoff-von Neumann decomposition, $\Pi_{\B}$ can be written as a convex combination of permutation matrices, i.e., there exist $c_1, \dots, c_d \in [0, 1]$ with $\sum_{p = 1}^{d} c_p = 1$ and $d$ different permutation matrices $P_{\B}^1, \dots, P_{\B}^d$ such that
\begin{align}\label{eq:proj_decom}
    \Pi_{\B} = \sum_{p = 1}^d c_p P_{\B}^p.
\end{align}
For any node $i \in [N]$, if $\Pi_{\B}^{ij} \neq 0$ then $j$ belongs to the same cluster of $i$, i.e., a particular cluster $\mathcal{K}_{\B}^l$. Since $\{\mathcal{K}_{\B}^l\}_{l = 1}^{k_{\B}}$ is an orbit partition, then there exists at least one permutation matrix $P_{\B}^p$ in the decomposition \eqref{eq:proj_decom} that permutes $i$ and $j$. Since the argument can be repeated for every node in the graph, we conclude that all permutation matrices $P_{\B}^p$ are elements of one of the equivalence classes $\mathcal{E}_{\B}^i$, and satisfy $P_{\B}^{p} L = L P_{\B}^{p}$. Then, $L \Pi_{\B} = \Pi_{\B} L$. 

Since commuting matrices preserve eigenspaces, then $L$ can be block-diagonalized using the set of eigenvectors of $\Pi_{\B}$. Repeating the same sequence of arguments to the adjacency matrix, the proposition holds. 
\end{proof}

\subsection{Stability characterization for the switching phenomenon}
\label{appendix:six_nodes}

Here we detail the stability characterization that underlies the switching phenomenon for an example of $6$ nodes in the top and bottom layers, see Figure~\ref{fig:directed_multilayer_n_12}. Figures~\ref{fig:directed_multilayer_n_12} a) and \ref{fig:directed_multilayer_n_12} b) display the different patterns as we vary $\alpha$ and $\beta$ while keeping $\sigma=0$. To demonstrate the validity of our numerical simulations, we also plot the largest of all transverse Lyapunov exponents ($\Lambda$) associated with each nontrivial cluster along with numerical simulations. A side-by-side comparison between synchronization errors $\e_{\mathcal{P}_{\T}}(t)$ and $\e_{\mathcal{P}_{\B}}(t)$ for different clusters and the corresponding $\Lambda_{\mathcal{K}_{\T}}$ and $\Lambda_{\mathcal{K}_{\B}}$ values consolidate the accuracy of numerical simulations. 

In Section~\ref{sec:bottom_layer_symm} we discussed the effect on invariant synchronous clusters in the bottom network when it is connected with the top: invariant clusters whose corresponding permutation matrices do not satisfy the conjugacy relation Equation~\eqref{eq:conjugacy_relation} are not invariant for any $\sigma>0$, see Section \ref{sec:restricting_symmetry}. Therefore, the top layer restricts the number of identically synchronized clusters, and consequently, pattern states at the bottom. For instance, in Figure~\ref{fig:directed_multilayer_n_12} b), when $\sigma=0$, the bottom layer exhibits patterns such as $\mathcal{P}_{\B}^2 = \{(a,b,b,c,c,c)\}$ and the complete synchronous pattern $\mathcal{P}_{\B}^3 = \{(a,a,a,a,a,a)\}$, confirmed by tracing the maximum Lyapunov exponents and the synchronization error $\mathfrak{e}$. However, these patterns are not invariant for any $\sigma>0$, or in other words, the corresponding permutation matrix $P_{\B}$ does not have a $P_{\T}$ in the top network that satisfies Equation~\eqref{eq:conjugacy_relation}. 

The top panel of Figure~\ref{fig:directed_multilayer_n_12} d) shows how the stability of bottom layers' patterns changes in the plane $(\alpha, \sigma)$ when driven by the top layer. The pattern $\mathcal{P}^0_{\mathfrak{B}}$ is stable at $\sigma=0$, and increasing $\alpha, \sigma$, the patterns $\mathcal{P}^4_{\mathfrak{B}}$ and $\mathcal{P}^5_{\mathfrak{B}}$ become stable. If we separate the isolated dynamics from coupling terms in Equation \eqref{eq:linear_flow_close}, the bottom panel of Figure~\ref{fig:directed_multilayer_n_12} d) shows the coupling matrix in the node space for the pattern state shown in Figure \ref{fig:directed_multilayer_n_12} c). To obtain this coupling matrix, we have assumed that $\alpha=\beta=\sigma=1$. The coupling matrix is such that $B^\perp$ is decoupled from $B^\parallel$ and $M^\perp$ from $M^\parallel$. However, due to the one-way \cite{lodi2021one} inter-layer dependence of the bottom network on top, $M^\perp$ is coupled with $B^\perp$, but the reverse is not true. Each row that is associated with $B^\perp$ and $M^\perp$ corresponds to a transverse perturbation that determines the stability of a cluster. Tracing all transverse Lyapunov exponents associated with a cluster with $\abs{\mathcal{K}^l_{(\T)\B}}>2$ provides a detailed analysis about which nodes in an invariant cluster can synchronize at a given $\alpha$ and $\beta$ \cite{Siddique2018}. For instance, for cluster $\mathcal{K}_{\T}^1$, the coupling matrix shows that if $x_i(0) \approx x_j(0)$ $\forall i,j \in \mathcal{K}_{\T}^1$, all nodes in this cluster synchronize simultaneously. 

There are two types of synchronized clusters: independent and intertwined. Independent clusters are those whose existence does not require any other cluster to exist, i.e., they can exist independently, whereas intertwined clusters require all clusters (that are intertwined) to exist simultaneously. If any cluster in a set of intertwined clusters desynchronizes, the stability of all intertwined clusters is lost as well. Using this fact, we destabilize a synchronized cluster in the bottom layer. Such independence or dependence of clusters is visible from the coupling matrix in Figure~\ref{fig:directed_multilayer_n_12} d) as well, which shows that clusters $\mathcal{K}_{\T}^1, \mathcal{K}_{\T}^2$ are independent, while $\mathcal{K}_{\B}^1, \mathcal{K}_{\B}^4$ are intertwined with $\mathcal{K}_{\T}^1, \mathcal{K}_{\T}^2$. While such intertwining of clusters is possible within the bottom layer also, our example contains intertwining across the layers only. 


\bibliography{aipsamp}

\begin{thebibliography}{58}%
\makeatletter
\providecommand \@ifxundefined [1]{%
 \@ifx{#1\undefined}
}%
\providecommand \@ifnum [1]{%
 \ifnum #1\expandafter \@firstoftwo
 \else \expandafter \@secondoftwo
 \fi
}%
\providecommand \@ifx [1]{%
 \ifx #1\expandafter \@firstoftwo
 \else \expandafter \@secondoftwo
 \fi
}%
\providecommand \natexlab [1]{#1}%
\providecommand \enquote  [1]{``#1''}%
\providecommand \bibnamefont  [1]{#1}%
\providecommand \bibfnamefont [1]{#1}%
\providecommand \citenamefont [1]{#1}%
\providecommand \href@noop [0]{\@secondoftwo}%
\providecommand \href [0]{\begingroup \@sanitize@url \@href}%
\providecommand \@href[1]{\@@startlink{#1}\@@href}%
\providecommand \@@href[1]{\endgroup#1\@@endlink}%
\providecommand \@sanitize@url [0]{\catcode `\\12\catcode `\$12\catcode `\&12\catcode `\#12\catcode `\^12\catcode `\_12\catcode `\%12\relax}%
\providecommand \@@startlink[1]{}%
\providecommand \@@endlink[0]{}%
\providecommand \url  [0]{\begingroup\@sanitize@url \@url }%
\providecommand \@url [1]{\endgroup\@href {#1}{\urlprefix }}%
\providecommand \urlprefix  [0]{URL }%
\providecommand \Eprint [0]{\href }%
\providecommand \doibase [0]{http://dx.doi.org/}%
\providecommand \selectlanguage [0]{\@gobble}%
\providecommand \bibinfo  [0]{\@secondoftwo}%
\providecommand \bibfield  [0]{\@secondoftwo}%
\providecommand \translation [1]{[#1]}%
\providecommand \BibitemOpen [0]{}%
\providecommand \bibitemStop [0]{}%
\providecommand \bibitemNoStop [0]{.\EOS\space}%
\providecommand \EOS [0]{\spacefactor3000\relax}%
\providecommand \BibitemShut  [1]{\csname bibitem#1\endcsname}%
\let\auto@bib@innerbib\@empty
\bibitem [{\citenamefont {Arenas}\ \emph {et~al.}(2008)\citenamefont {Arenas}, \citenamefont {Díaz-Guilera}, \citenamefont {Kurths}, \citenamefont {Moreno},\ and\ \citenamefont {Zhou}}]{Arenas2008}%
  \BibitemOpen
  \bibfield  {author} {\bibinfo {author} {\bibfnamefont {A.}~\bibnamefont {Arenas}}, \bibinfo {author} {\bibfnamefont {A.}~\bibnamefont {Díaz-Guilera}}, \bibinfo {author} {\bibfnamefont {J.}~\bibnamefont {Kurths}}, \bibinfo {author} {\bibfnamefont {Y.}~\bibnamefont {Moreno}}, \ and\ \bibinfo {author} {\bibfnamefont {C.}~\bibnamefont {Zhou}},\ }\bibfield  {title} {\enquote {\bibinfo {title} {Synchronization in complex networks},}\ }\href {\doibase https://doi.org/10.1016/j.physrep.2008.09.002} {\bibfield  {journal} {\bibinfo  {journal} {Physics Reports}\ }\textbf {\bibinfo {volume} {469}},\ \bibinfo {pages} {93--153} (\bibinfo {year} {2008})}\BibitemShut {NoStop}%
\bibitem [{\citenamefont {Rodrigues}\ \emph {et~al.}(2016)\citenamefont {Rodrigues}, \citenamefont {Peron}, \citenamefont {Ji},\ and\ \citenamefont {Kurths}}]{Rodrigues2016}%
  \BibitemOpen
  \bibfield  {author} {\bibinfo {author} {\bibfnamefont {F.~A.}\ \bibnamefont {Rodrigues}}, \bibinfo {author} {\bibfnamefont {T.~K.~D.}\ \bibnamefont {Peron}}, \bibinfo {author} {\bibfnamefont {P.}~\bibnamefont {Ji}}, \ and\ \bibinfo {author} {\bibfnamefont {J.}~\bibnamefont {Kurths}},\ }\bibfield  {title} {\enquote {\bibinfo {title} {{The Kuramoto model in complex networks}},}\ }\href {\doibase https://doi.org/10.1016/j.physrep.2015.10.008} {\bibfield  {journal} {\bibinfo  {journal} {Physics Reports}\ }\textbf {\bibinfo {volume} {610}},\ \bibinfo {pages} {1--98} (\bibinfo {year} {2016})},\ \bibinfo {note} {the Kuramoto model in complex networks}\BibitemShut {NoStop}%
\bibitem [{\citenamefont {Panaggio}\ and\ \citenamefont {Abrams}(2015)}]{Panaggio_2015}%
  \BibitemOpen
  \bibfield  {author} {\bibinfo {author} {\bibfnamefont {M.~J.}\ \bibnamefont {Panaggio}}\ and\ \bibinfo {author} {\bibfnamefont {D.~M.}\ \bibnamefont {Abrams}},\ }\bibfield  {title} {\enquote {\bibinfo {title} {Chimera states: coexistence of coherence and incoherence in networks of coupled oscillators},}\ }\href {\doibase 10.1088/0951-7715/28/3/R67} {\bibfield  {journal} {\bibinfo  {journal} {Nonlinearity}\ }\textbf {\bibinfo {volume} {28}},\ \bibinfo {pages} {R67} (\bibinfo {year} {2015})}\BibitemShut {NoStop}%
\bibitem [{\citenamefont {Wang}\ and\ \citenamefont {Liu}(2020)}]{Wang2020}%
  \BibitemOpen
  \bibfield  {author} {\bibinfo {author} {\bibfnamefont {Z.}~\bibnamefont {Wang}}\ and\ \bibinfo {author} {\bibfnamefont {Z.}~\bibnamefont {Liu}},\ }\bibfield  {title} {\enquote {\bibinfo {title} {A brief review of chimera state in empirical brain networks},}\ }\href {\doibase 10.3389/fphys.2020.00724} {\bibfield  {journal} {\bibinfo  {journal} {Frontiers in Physiology}\ }\textbf {\bibinfo {volume} {11}} (\bibinfo {year} {2020}),\ 10.3389/fphys.2020.00724}\BibitemShut {NoStop}%
\bibitem [{\citenamefont {Pecora}\ \emph {et~al.}(2014)\citenamefont {Pecora}, \citenamefont {Sorrentino}, \citenamefont {Hagerstrom}, \citenamefont {Murphy},\ and\ \citenamefont {Roy}}]{Pecora2014}%
  \BibitemOpen
  \bibfield  {author} {\bibinfo {author} {\bibfnamefont {L.~M.}\ \bibnamefont {Pecora}}, \bibinfo {author} {\bibfnamefont {F.}~\bibnamefont {Sorrentino}}, \bibinfo {author} {\bibfnamefont {A.~M.}\ \bibnamefont {Hagerstrom}}, \bibinfo {author} {\bibfnamefont {T.~E.}\ \bibnamefont {Murphy}}, \ and\ \bibinfo {author} {\bibfnamefont {R.}~\bibnamefont {Roy}},\ }\bibfield  {title} {\enquote {\bibinfo {title} {Cluster synchronization and isolated desynchronization in complex networks with symmetries},}\ }\href {\doibase 10.1038/ncomms5079} {\bibfield  {journal} {\bibinfo  {journal} {Nature Communications}\ }\textbf {\bibinfo {volume} {5}},\ \bibinfo {pages} {4079} (\bibinfo {year} {2014})}\BibitemShut {NoStop}%
\bibitem [{\citenamefont {Parastesh}\ \emph {et~al.}(2021)\citenamefont {Parastesh}, \citenamefont {Jafari}, \citenamefont {Azarnoush}, \citenamefont {Shahriari}, \citenamefont {Wang}, \citenamefont {Boccaletti},\ and\ \citenamefont {Perc}}]{Parastesh2021}%
  \BibitemOpen
  \bibfield  {author} {\bibinfo {author} {\bibfnamefont {F.}~\bibnamefont {Parastesh}}, \bibinfo {author} {\bibfnamefont {S.}~\bibnamefont {Jafari}}, \bibinfo {author} {\bibfnamefont {H.}~\bibnamefont {Azarnoush}}, \bibinfo {author} {\bibfnamefont {Z.}~\bibnamefont {Shahriari}}, \bibinfo {author} {\bibfnamefont {Z.}~\bibnamefont {Wang}}, \bibinfo {author} {\bibfnamefont {S.}~\bibnamefont {Boccaletti}}, \ and\ \bibinfo {author} {\bibfnamefont {M.}~\bibnamefont {Perc}},\ }\bibfield  {title} {\enquote {\bibinfo {title} {Chimeras},}\ }\href {\doibase https://doi.org/10.1016/j.physrep.2020.10.003} {\bibfield  {journal} {\bibinfo  {journal} {Physics Reports}\ }\textbf {\bibinfo {volume} {898}},\ \bibinfo {pages} {1--114} (\bibinfo {year} {2021})},\ \bibinfo {note} {chimeras}\BibitemShut {NoStop}%
\bibitem [{\citenamefont {Buzsaki}(2006)}]{buzsaki2006rhythms}%
  \BibitemOpen
  \bibfield  {author} {\bibinfo {author} {\bibfnamefont {G.}~\bibnamefont {Buzsaki}},\ }\href@noop {} {\emph {\bibinfo {title} {Rhythms of the Brain}}}\ (\bibinfo  {publisher} {Oxford university press},\ \bibinfo {year} {2006})\BibitemShut {NoStop}%
\bibitem [{\citenamefont {Tognoli}\ and\ \citenamefont {Kelso}(2014)}]{Tognoli2014}%
  \BibitemOpen
  \bibfield  {author} {\bibinfo {author} {\bibfnamefont {E.}~\bibnamefont {Tognoli}}\ and\ \bibinfo {author} {\bibfnamefont {J.}~\bibnamefont {Kelso}},\ }\bibfield  {title} {\enquote {\bibinfo {title} {The metastable brain},}\ }\href {\doibase 10.1016/j.neuron.2013.12.022} {\bibfield  {journal} {\bibinfo  {journal} {Neuron}\ }\textbf {\bibinfo {volume} {81}},\ \bibinfo {pages} {35--48} (\bibinfo {year} {2014})}\BibitemShut {NoStop}%
\bibitem [{\citenamefont {Deco}\ \emph {et~al.}(2015)\citenamefont {Deco}, \citenamefont {Tononi}, \citenamefont {Boly},\ and\ \citenamefont {Kringelbach}}]{Deco2015}%
  \BibitemOpen
  \bibfield  {author} {\bibinfo {author} {\bibfnamefont {G.}~\bibnamefont {Deco}}, \bibinfo {author} {\bibfnamefont {G.}~\bibnamefont {Tononi}}, \bibinfo {author} {\bibfnamefont {M.}~\bibnamefont {Boly}}, \ and\ \bibinfo {author} {\bibfnamefont {M.~L.}\ \bibnamefont {Kringelbach}},\ }\bibfield  {title} {\enquote {\bibinfo {title} {Rethinking segregation and integration: contributions of whole-brain modelling},}\ }\href {\doibase 10.1038/nrn3963} {\bibfield  {journal} {\bibinfo  {journal} {Nature Reviews Neuroscience}\ }\textbf {\bibinfo {volume} {16}},\ \bibinfo {pages} {430--439} (\bibinfo {year} {2015})}\BibitemShut {NoStop}%
\bibitem [{\citenamefont {Handwerker}\ \emph {et~al.}(2012)\citenamefont {Handwerker}, \citenamefont {Roopchansingh}, \citenamefont {Gonzalez-Castillo},\ and\ \citenamefont {Bandettini}}]{handwerker2012}%
  \BibitemOpen
  \bibfield  {author} {\bibinfo {author} {\bibfnamefont {D.~A.}\ \bibnamefont {Handwerker}}, \bibinfo {author} {\bibfnamefont {V.}~\bibnamefont {Roopchansingh}}, \bibinfo {author} {\bibfnamefont {J.}~\bibnamefont {Gonzalez-Castillo}}, \ and\ \bibinfo {author} {\bibfnamefont {P.~A.}\ \bibnamefont {Bandettini}},\ }\bibfield  {title} {\enquote {\bibinfo {title} {Periodic changes in fmri connectivity},}\ }\href@noop {} {\bibfield  {journal} {\bibinfo  {journal} {Neuroimage}\ }\textbf {\bibinfo {volume} {63}},\ \bibinfo {pages} {1712--1719} (\bibinfo {year} {2012})}\BibitemShut {NoStop}%
\bibitem [{\citenamefont {Chang}\ and\ \citenamefont {Glover}(2010)}]{chang2010}%
  \BibitemOpen
  \bibfield  {author} {\bibinfo {author} {\bibfnamefont {C.}~\bibnamefont {Chang}}\ and\ \bibinfo {author} {\bibfnamefont {G.~H.}\ \bibnamefont {Glover}},\ }\bibfield  {title} {\enquote {\bibinfo {title} {Time--frequency dynamics of resting-state brain connectivity measured with fmri},}\ }\href@noop {} {\bibfield  {journal} {\bibinfo  {journal} {Neuroimage}\ }\textbf {\bibinfo {volume} {50}},\ \bibinfo {pages} {81--98} (\bibinfo {year} {2010})}\BibitemShut {NoStop}%
\bibitem [{\citenamefont {Allen}\ \emph {et~al.}(2014)\citenamefont {Allen}, \citenamefont {Damaraju}, \citenamefont {Plis}, \citenamefont {Erhardt}, \citenamefont {Eichele},\ and\ \citenamefont {Calhoun}}]{allen2014}%
  \BibitemOpen
  \bibfield  {author} {\bibinfo {author} {\bibfnamefont {E.~A.}\ \bibnamefont {Allen}}, \bibinfo {author} {\bibfnamefont {E.}~\bibnamefont {Damaraju}}, \bibinfo {author} {\bibfnamefont {S.~M.}\ \bibnamefont {Plis}}, \bibinfo {author} {\bibfnamefont {E.~B.}\ \bibnamefont {Erhardt}}, \bibinfo {author} {\bibfnamefont {T.}~\bibnamefont {Eichele}}, \ and\ \bibinfo {author} {\bibfnamefont {V.~D.}\ \bibnamefont {Calhoun}},\ }\bibfield  {title} {\enquote {\bibinfo {title} {Tracking whole-brain connectivity dynamics in the resting state},}\ }\href@noop {} {\bibfield  {journal} {\bibinfo  {journal} {Cerebral cortex}\ }\textbf {\bibinfo {volume} {24}},\ \bibinfo {pages} {663--676} (\bibinfo {year} {2014})}\BibitemShut {NoStop}%
\bibitem [{\citenamefont {Stewart}, \citenamefont {Golubitsky},\ and\ \citenamefont {Pivato}(2003)}]{Stewart_2003}%
  \BibitemOpen
  \bibfield  {author} {\bibinfo {author} {\bibfnamefont {I.}~\bibnamefont {Stewart}}, \bibinfo {author} {\bibfnamefont {M.}~\bibnamefont {Golubitsky}}, \ and\ \bibinfo {author} {\bibfnamefont {M.}~\bibnamefont {Pivato}},\ }\bibfield  {title} {\enquote {\bibinfo {title} {Symmetry groupoids and patterns of synchrony in coupled cell networks},}\ }\href {\doibase 10.1137/S1111111103419896} {\bibfield  {journal} {\bibinfo  {journal} {SIAM Journal on Applied Dynamical Systems}\ }\textbf {\bibinfo {volume} {2}},\ \bibinfo {pages} {609--646} (\bibinfo {year} {2003})},\ \Eprint {http://arxiv.org/abs/https://doi.org/10.1137/S1111111103419896} {https://doi.org/10.1137/S1111111103419896} \BibitemShut {NoStop}%
\bibitem [{\citenamefont {Golubitsky}, \citenamefont {Stewart},\ and\ \citenamefont {T\"{o}r\"{o}k}(2005)}]{Golubitsky2005}%
  \BibitemOpen
  \bibfield  {author} {\bibinfo {author} {\bibfnamefont {M.}~\bibnamefont {Golubitsky}}, \bibinfo {author} {\bibfnamefont {I.}~\bibnamefont {Stewart}}, \ and\ \bibinfo {author} {\bibfnamefont {A.}~\bibnamefont {T\"{o}r\"{o}k}},\ }\bibfield  {title} {\enquote {\bibinfo {title} {Patterns of synchrony in coupled cell networks with multiple arrows},}\ }\href {\doibase 10.1137/040612634} {\bibfield  {journal} {\bibinfo  {journal} {SIAM Journal on Applied Dynamical Systems}\ }\textbf {\bibinfo {volume} {4}},\ \bibinfo {pages} {78--100} (\bibinfo {year} {2005})},\ \Eprint {http://arxiv.org/abs/https://doi.org/10.1137/040612634} {https://doi.org/10.1137/040612634} \BibitemShut {NoStop}%
\bibitem [{\citenamefont {Belykh}\ \emph {et~al.}(2008)\citenamefont {Belykh}, \citenamefont {Osipov}, \citenamefont {Petrov}, \citenamefont {Suykens},\ and\ \citenamefont {Vandewalle}}]{Belykh_2008}%
  \BibitemOpen
  \bibfield  {author} {\bibinfo {author} {\bibfnamefont {V.~N.}\ \bibnamefont {Belykh}}, \bibinfo {author} {\bibfnamefont {G.~V.}\ \bibnamefont {Osipov}}, \bibinfo {author} {\bibfnamefont {V.~S.}\ \bibnamefont {Petrov}}, \bibinfo {author} {\bibfnamefont {J.~A.~K.}\ \bibnamefont {Suykens}}, \ and\ \bibinfo {author} {\bibfnamefont {J.}~\bibnamefont {Vandewalle}},\ }\bibfield  {title} {\enquote {\bibinfo {title} {{Cluster synchronization in oscillatory networks}},}\ }\href {\doibase 10.1063/1.2956986} {\bibfield  {journal} {\bibinfo  {journal} {Chaos: An Interdisciplinary Journal of Nonlinear Science}\ }\textbf {\bibinfo {volume} {18}},\ \bibinfo {pages} {037106} (\bibinfo {year} {2008})},\ \Eprint {http://arxiv.org/abs/https://pubs.aip.org/aip/cha/article-pdf/doi/10.1063/1.2956986/13857117/037106\_1\_online.pdf} {https://pubs.aip.org/aip/cha/article-pdf/doi/10.1063/1.2956986/13857117/037106\_1\_online.pdf} \BibitemShut {NoStop}%
\bibitem [{\citenamefont {Golubitsky}(2002)}]{GolubitskyMartin2002Tspf}%
  \BibitemOpen
  \bibfield  {author} {\bibinfo {author} {\bibfnamefont {M.}~\bibnamefont {Golubitsky}},\ }\href@noop {} {\emph {\bibinfo {title} {The symmetry perspective from equilibrium to chaos in phase space and physical space}}},\ Progress in mathematics (Boston, Mass.) v. 200\ (\bibinfo  {publisher} {Birkhäuser},\ \bibinfo {address} {Basel Boston},\ \bibinfo {year} {2002})\BibitemShut {NoStop}%
\bibitem [{\citenamefont {Sorrentino}\ \emph {et~al.}(2016)\citenamefont {Sorrentino}, \citenamefont {Pecora}, \citenamefont {Hagerstrom}, \citenamefont {Murphy},\ and\ \citenamefont {Roy}}]{Sorrentino_2016}%
  \BibitemOpen
  \bibfield  {author} {\bibinfo {author} {\bibfnamefont {F.}~\bibnamefont {Sorrentino}}, \bibinfo {author} {\bibfnamefont {L.~M.}\ \bibnamefont {Pecora}}, \bibinfo {author} {\bibfnamefont {A.~M.}\ \bibnamefont {Hagerstrom}}, \bibinfo {author} {\bibfnamefont {T.~E.}\ \bibnamefont {Murphy}}, \ and\ \bibinfo {author} {\bibfnamefont {R.}~\bibnamefont {Roy}},\ }\bibfield  {title} {\enquote {\bibinfo {title} {Complete characterization of the stability of cluster synchronization in complex dynamical networks},}\ }\href {\doibase 10.1126/sciadv.1501737} {\bibfield  {journal} {\bibinfo  {journal} {Science Advances}\ }\textbf {\bibinfo {volume} {2}},\ \bibinfo {pages} {e1501737} (\bibinfo {year} {2016})},\ \Eprint {http://arxiv.org/abs/https://www.science.org/doi/pdf/10.1126/sciadv.1501737} {https://www.science.org/doi/pdf/10.1126/sciadv.1501737} \BibitemShut {NoStop}%
\bibitem [{\citenamefont {Koch}\ and\ \citenamefont {Laurent}(1999)}]{Koch_1999}%
  \BibitemOpen
  \bibfield  {author} {\bibinfo {author} {\bibfnamefont {C.}~\bibnamefont {Koch}}\ and\ \bibinfo {author} {\bibfnamefont {G.}~\bibnamefont {Laurent}},\ }\bibfield  {title} {\enquote {\bibinfo {title} {Complexity and the nervous system},}\ }\href {\doibase 10.1126/science.284.5411.96} {\bibfield  {journal} {\bibinfo  {journal} {Science}\ }\textbf {\bibinfo {volume} {284}},\ \bibinfo {pages} {96--98} (\bibinfo {year} {1999})},\ \Eprint {http://arxiv.org/abs/https://www.science.org/doi/pdf/10.1126/science.284.5411.96} {https://www.science.org/doi/pdf/10.1126/science.284.5411.96} \BibitemShut {NoStop}%
\bibitem [{\citenamefont {Bullmore}\ and\ \citenamefont {Sporns}(2009)}]{bullmore2009}%
  \BibitemOpen
  \bibfield  {author} {\bibinfo {author} {\bibfnamefont {E.}~\bibnamefont {Bullmore}}\ and\ \bibinfo {author} {\bibfnamefont {O.}~\bibnamefont {Sporns}},\ }\bibfield  {title} {\enquote {\bibinfo {title} {Complex brain networks: graph theoretical analysis of structural and functional systems},}\ }\href@noop {} {\bibfield  {journal} {\bibinfo  {journal} {Nature Reviews Neuroscience}\ }\textbf {\bibinfo {volume} {10}},\ \bibinfo {pages} {186} (\bibinfo {year} {2009})}\BibitemShut {NoStop}%
\bibitem [{\citenamefont {De~Domenico}\ \emph {et~al.}(2013)\citenamefont {De~Domenico}, \citenamefont {Sol\'e-Ribalta}, \citenamefont {Cozzo}, \citenamefont {Kivel\"a}, \citenamefont {Moreno}, \citenamefont {Porter}, \citenamefont {G\'omez},\ and\ \citenamefont {Arenas}}]{Domenico2013}%
  \BibitemOpen
  \bibfield  {author} {\bibinfo {author} {\bibfnamefont {M.}~\bibnamefont {De~Domenico}}, \bibinfo {author} {\bibfnamefont {A.}~\bibnamefont {Sol\'e-Ribalta}}, \bibinfo {author} {\bibfnamefont {E.}~\bibnamefont {Cozzo}}, \bibinfo {author} {\bibfnamefont {M.}~\bibnamefont {Kivel\"a}}, \bibinfo {author} {\bibfnamefont {Y.}~\bibnamefont {Moreno}}, \bibinfo {author} {\bibfnamefont {M.~A.}\ \bibnamefont {Porter}}, \bibinfo {author} {\bibfnamefont {S.}~\bibnamefont {G\'omez}}, \ and\ \bibinfo {author} {\bibfnamefont {A.}~\bibnamefont {Arenas}},\ }\bibfield  {title} {\enquote {\bibinfo {title} {Mathematical formulation of multilayer networks},}\ }\href {\doibase 10.1103/PhysRevX.3.041022} {\bibfield  {journal} {\bibinfo  {journal} {Phys. Rev. X}\ }\textbf {\bibinfo {volume} {3}},\ \bibinfo {pages} {041022} (\bibinfo {year} {2013})}\BibitemShut {NoStop}%
\bibitem [{\citenamefont {Crofts}, \citenamefont {Forrester},\ and\ \citenamefont {O'Dea}(2016)}]{Crofts_2016}%
  \BibitemOpen
  \bibfield  {author} {\bibinfo {author} {\bibfnamefont {J.~J.}\ \bibnamefont {Crofts}}, \bibinfo {author} {\bibfnamefont {M.}~\bibnamefont {Forrester}}, \ and\ \bibinfo {author} {\bibfnamefont {R.~D.}\ \bibnamefont {O'Dea}},\ }\bibfield  {title} {\enquote {\bibinfo {title} {Structure-function clustering in multiplex brain networks},}\ }\href {\doibase 10.1209/0295-5075/116/18003} {\bibfield  {journal} {\bibinfo  {journal} {Europhysics Letters}\ }\textbf {\bibinfo {volume} {116}},\ \bibinfo {pages} {18003} (\bibinfo {year} {2016})}\BibitemShut {NoStop}%
\bibitem [{\citenamefont {De~Domenico}(2017)}]{DeDomenico_2017}%
  \BibitemOpen
  \bibfield  {author} {\bibinfo {author} {\bibfnamefont {M.}~\bibnamefont {De~Domenico}},\ }\bibfield  {title} {\enquote {\bibinfo {title} {{Multilayer modeling and analysis of human brain networks}},}\ }\href {\doibase 10.1093/gigascience/gix004} {\bibfield  {journal} {\bibinfo  {journal} {GigaScience}\ }\textbf {\bibinfo {volume} {6}},\ \bibinfo {pages} {gix004} (\bibinfo {year} {2017})},\ \Eprint {http://arxiv.org/abs/https://academic.oup.com/gigascience/article-pdf/6/5/gix004/25514568/gix004\_reviewer\_3\_report\_(original\_submission).pdf} {https://academic.oup.com/gigascience/article-pdf/6/5/gix004/25514568/gix004\_reviewer\_3\_report\_(original\_submission).pdf} \BibitemShut {NoStop}%
\bibitem [{\citenamefont {Battiston}\ \emph {et~al.}(2017)\citenamefont {Battiston}, \citenamefont {Nicosia}, \citenamefont {Chavez},\ and\ \citenamefont {Latora}}]{Battiston_2017}%
  \BibitemOpen
  \bibfield  {author} {\bibinfo {author} {\bibfnamefont {F.}~\bibnamefont {Battiston}}, \bibinfo {author} {\bibfnamefont {V.}~\bibnamefont {Nicosia}}, \bibinfo {author} {\bibfnamefont {M.}~\bibnamefont {Chavez}}, \ and\ \bibinfo {author} {\bibfnamefont {V.}~\bibnamefont {Latora}},\ }\bibfield  {title} {\enquote {\bibinfo {title} {{Multilayer motif analysis of brain networks}},}\ }\href {\doibase 10.1063/1.4979282} {\bibfield  {journal} {\bibinfo  {journal} {Chaos: An Interdisciplinary Journal of Nonlinear Science}\ }\textbf {\bibinfo {volume} {27}},\ \bibinfo {pages} {047404} (\bibinfo {year} {2017})},\ \Eprint {http://arxiv.org/abs/https://pubs.aip.org/aip/cha/article-pdf/doi/10.1063/1.4979282/13299057/047404\_1\_online.pdf} {https://pubs.aip.org/aip/cha/article-pdf/doi/10.1063/1.4979282/13299057/047404\_1\_online.pdf} \BibitemShut {NoStop}%
\bibitem [{\citenamefont {Frolov}, \citenamefont {Maksimenko},\ and\ \citenamefont {Hramov}(2020)}]{Frolov_2020}%
  \BibitemOpen
  \bibfield  {author} {\bibinfo {author} {\bibfnamefont {N.}~\bibnamefont {Frolov}}, \bibinfo {author} {\bibfnamefont {V.}~\bibnamefont {Maksimenko}}, \ and\ \bibinfo {author} {\bibfnamefont {A.}~\bibnamefont {Hramov}},\ }\bibfield  {title} {\enquote {\bibinfo {title} {{Revealing a multiplex brain network through the analysis of recurrences}},}\ }\href {\doibase 10.1063/5.0028053} {\bibfield  {journal} {\bibinfo  {journal} {Chaos: An Interdisciplinary Journal of Nonlinear Science}\ }\textbf {\bibinfo {volume} {30}},\ \bibinfo {pages} {121108} (\bibinfo {year} {2020})},\ \Eprint {http://arxiv.org/abs/https://pubs.aip.org/aip/cha/article-pdf/doi/10.1063/5.0028053/14106595/121108\_1\_online.pdf} {https://pubs.aip.org/aip/cha/article-pdf/doi/10.1063/5.0028053/14106595/121108\_1\_online.pdf} \BibitemShut {NoStop}%
\bibitem [{\citenamefont {Vaiana}\ and\ \citenamefont {Muldoon}(2020)}]{Vaiana2020}%
  \BibitemOpen
  \bibfield  {author} {\bibinfo {author} {\bibfnamefont {M.}~\bibnamefont {Vaiana}}\ and\ \bibinfo {author} {\bibfnamefont {S.~F.}\ \bibnamefont {Muldoon}},\ }\bibfield  {title} {\enquote {\bibinfo {title} {Multilayer brain networks},}\ }\href {\doibase 10.1007/s00332-017-9436-8} {\bibfield  {journal} {\bibinfo  {journal} {Journal of Nonlinear Science}\ }\textbf {\bibinfo {volume} {30}},\ \bibinfo {pages} {2147--2169} (\bibinfo {year} {2020})}\BibitemShut {NoStop}%
\bibitem [{\citenamefont {Della~Rossa}\ \emph {et~al.}(2020)\citenamefont {Della~Rossa}, \citenamefont {Pecora}, \citenamefont {Blaha}, \citenamefont {Shirin}, \citenamefont {Klickstein},\ and\ \citenamefont {Sorrentino}}]{DellaRossa2020}%
  \BibitemOpen
  \bibfield  {author} {\bibinfo {author} {\bibfnamefont {F.}~\bibnamefont {Della~Rossa}}, \bibinfo {author} {\bibfnamefont {L.}~\bibnamefont {Pecora}}, \bibinfo {author} {\bibfnamefont {K.}~\bibnamefont {Blaha}}, \bibinfo {author} {\bibfnamefont {A.}~\bibnamefont {Shirin}}, \bibinfo {author} {\bibfnamefont {I.}~\bibnamefont {Klickstein}}, \ and\ \bibinfo {author} {\bibfnamefont {F.}~\bibnamefont {Sorrentino}},\ }\bibfield  {title} {\enquote {\bibinfo {title} {Symmetries and cluster synchronization in multilayer networks},}\ }\href {\doibase 10.1038/s41467-020-16343-0} {\bibfield  {journal} {\bibinfo  {journal} {Nature Communications}\ }\textbf {\bibinfo {volume} {11}},\ \bibinfo {pages} {3179} (\bibinfo {year} {2020})}\BibitemShut {NoStop}%
\bibitem [{\citenamefont {BELYKH}, \citenamefont {BELYKH},\ and\ \citenamefont {MOSEKILDE}(2005)}]{BELYKH_2005}%
  \BibitemOpen
  \bibfield  {author} {\bibinfo {author} {\bibfnamefont {V.}~\bibnamefont {BELYKH}}, \bibinfo {author} {\bibfnamefont {I.}~\bibnamefont {BELYKH}}, \ and\ \bibinfo {author} {\bibfnamefont {E.}~\bibnamefont {MOSEKILDE}},\ }\bibfield  {title} {\enquote {\bibinfo {title} {Hyperbolic plykin attractor can exist in neuron models},}\ }\href {\doibase 10.1142/S0218127405014222} {\bibfield  {journal} {\bibinfo  {journal} {International Journal of Bifurcation and Chaos}\ }\textbf {\bibinfo {volume} {15}},\ \bibinfo {pages} {3567--3578} (\bibinfo {year} {2005})},\ \Eprint {http://arxiv.org/abs/https://doi.org/10.1142/S0218127405014222} {https://doi.org/10.1142/S0218127405014222} \BibitemShut {NoStop}%
\bibitem [{\citenamefont {SHILNIKOV}\ and\ \citenamefont {KOLOMIETS}(2008)}]{SHILNIKOV_2008}%
  \BibitemOpen
  \bibfield  {author} {\bibinfo {author} {\bibfnamefont {A.}~\bibnamefont {SHILNIKOV}}\ and\ \bibinfo {author} {\bibfnamefont {M.}~\bibnamefont {KOLOMIETS}},\ }\bibfield  {title} {\enquote {\bibinfo {title} {Methods of the qualitative theory for the hindmarsh–rose model: A case study – a tutorial},}\ }\href {\doibase 10.1142/S0218127408021634} {\bibfield  {journal} {\bibinfo  {journal} {International Journal of Bifurcation and Chaos}\ }\textbf {\bibinfo {volume} {18}},\ \bibinfo {pages} {2141--2168} (\bibinfo {year} {2008})},\ \Eprint {http://arxiv.org/abs/https://doi.org/10.1142/S0218127408021634} {https://doi.org/10.1142/S0218127408021634} \BibitemShut {NoStop}%
\bibitem [{\citenamefont {Harary}(1969)}]{HararyFrank1969Gt}%
  \BibitemOpen
  \bibfield  {author} {\bibinfo {author} {\bibfnamefont {F.}~\bibnamefont {Harary}},\ }\href@noop {} {\emph {\bibinfo {title} {Graph theory}}},\ Addison-Wesley series in mathematics\ (\bibinfo  {publisher} {Addison-Wesley Pub. Co},\ \bibinfo {address} {Reading, Mass},\ \bibinfo {year} {1969})\BibitemShut {NoStop}%
\bibitem [{\citenamefont {Schaub}\ \emph {et~al.}(2016)\citenamefont {Schaub}, \citenamefont {O'Clery}, \citenamefont {Billeh}, \citenamefont {Delvenne}, \citenamefont {Lambiotte},\ and\ \citenamefont {Barahona}}]{Schaub2016}%
  \BibitemOpen
  \bibfield  {author} {\bibinfo {author} {\bibfnamefont {M.~T.}\ \bibnamefont {Schaub}}, \bibinfo {author} {\bibfnamefont {N.}~\bibnamefont {O'Clery}}, \bibinfo {author} {\bibfnamefont {Y.~N.}\ \bibnamefont {Billeh}}, \bibinfo {author} {\bibfnamefont {J.-C.}\ \bibnamefont {Delvenne}}, \bibinfo {author} {\bibfnamefont {R.}~\bibnamefont {Lambiotte}}, \ and\ \bibinfo {author} {\bibfnamefont {M.}~\bibnamefont {Barahona}},\ }\bibfield  {title} {\enquote {\bibinfo {title} {{Graph partitions and cluster synchronization in networks of oscillators}},}\ }\href {\doibase 10.1063/1.4961065} {\bibfield  {journal} {\bibinfo  {journal} {Chaos: An Interdisciplinary Journal of Nonlinear Science}\ }\textbf {\bibinfo {volume} {26}},\ \bibinfo {pages} {094821} (\bibinfo {year} {2016})},\ \Eprint {http://arxiv.org/abs/https://pubs.aip.org/aip/cha/article-pdf/doi/10.1063/1.4961065/14614631/094821\_1\_online.pdf} {https://pubs.aip.org/aip/cha/article-pdf/doi/10.1063/1.4961065/14614631/094821\_1\_online.pdf} \BibitemShut {NoStop}%
\bibitem [{\citenamefont {DeVille}\ and\ \citenamefont {Lerman}(2015)}]{DeVille2015}%
  \BibitemOpen
  \bibfield  {author} {\bibinfo {author} {\bibfnamefont {L.}~\bibnamefont {DeVille}}\ and\ \bibinfo {author} {\bibfnamefont {E.}~\bibnamefont {Lerman}},\ }\bibfield  {title} {\enquote {\bibinfo {title} {Modular dynamical systems on networks},}\ }\href {http://eudml.org/doc/277581} {\bibfield  {journal} {\bibinfo  {journal} {Journal of the European Mathematical Society}\ }\textbf {\bibinfo {volume} {017}},\ \bibinfo {pages} {2977--3013} (\bibinfo {year} {2015})}\BibitemShut {NoStop}%
\bibitem [{\citenamefont {Kudose}()}]{orbit_partition}%
  \BibitemOpen
  \bibfield  {author} {\bibinfo {author} {\bibfnamefont {S.}~\bibnamefont {Kudose}},\ }\href@noop {} {\enquote {\bibinfo {title} {{Equitable partitions and orbit partitions}},}\ }\bibinfo {howpublished} {https://www.math.uchicago.edu/~may/VIGRE/VIGRE2009/REUPapers/Kudose.pdf},\ \bibinfo {note} {unpublished}\BibitemShut {NoStop}%
\bibitem [{\citenamefont {Stein}\ and\ \citenamefont {Joyner}(2005)}]{sage}%
  \BibitemOpen
  \bibfield  {author} {\bibinfo {author} {\bibfnamefont {W.}~\bibnamefont {Stein}}\ and\ \bibinfo {author} {\bibfnamefont {D.}~\bibnamefont {Joyner}},\ }\href@noop {} {\enquote {\bibinfo {title} {{SAGE: Software for Algebra and Geometry Experimentation}},}\ }\bibinfo {howpublished} {https://www.sagemath.org} (\bibinfo {year} {2005})\BibitemShut {NoStop}%
\bibitem [{\citenamefont {Lodi}, \citenamefont {Sorrentino},\ and\ \citenamefont {Storace}(2021)}]{lodi2021one}%
  \BibitemOpen
  \bibfield  {author} {\bibinfo {author} {\bibfnamefont {M.}~\bibnamefont {Lodi}}, \bibinfo {author} {\bibfnamefont {F.}~\bibnamefont {Sorrentino}}, \ and\ \bibinfo {author} {\bibfnamefont {M.}~\bibnamefont {Storace}},\ }\bibfield  {title} {\enquote {\bibinfo {title} {One-way dependent clusters and stability of cluster synchronization in directed networks},}\ }\href@noop {} {\bibfield  {journal} {\bibinfo  {journal} {Nature communications}\ }\textbf {\bibinfo {volume} {12}},\ \bibinfo {pages} {4073} (\bibinfo {year} {2021})}\BibitemShut {NoStop}%
\bibitem [{\citenamefont {Hahn}\ and\ \citenamefont {Sabidussi}(1997)}]{HahnGena1997GSAM}%
  \BibitemOpen
  \bibfield  {author} {\bibinfo {author} {\bibfnamefont {G.}~\bibnamefont {Hahn}}\ and\ \bibinfo {author} {\bibfnamefont {G.}~\bibnamefont {Sabidussi}},\ }\href@noop {} {\emph {\bibinfo {title} {Graph Symmetry}}},\ \bibinfo {series} {Nato Science Series C:, Mathematical and Physical Sciences}, Vol.\ \bibinfo {volume} {497}\ (\bibinfo  {publisher} {Springer Netherlands},\ \bibinfo {address} {Dordrecht},\ \bibinfo {year} {1997})\BibitemShut {NoStop}%
\bibitem [{Note1()}]{Note1}%
  \BibitemOpen
  \bibinfo {note} {The network permutation matrix in \cite {Lin_2016}.}\BibitemShut {Stop}%
\bibitem [{\citenamefont {Pereira}\ \emph {et~al.}(2014)\citenamefont {Pereira}, \citenamefont {Eldering}, \citenamefont {Rasmussen},\ and\ \citenamefont {Veneziani}}]{Pereira_2014}%
  \BibitemOpen
  \bibfield  {author} {\bibinfo {author} {\bibfnamefont {T.}~\bibnamefont {Pereira}}, \bibinfo {author} {\bibfnamefont {J.}~\bibnamefont {Eldering}}, \bibinfo {author} {\bibfnamefont {M.}~\bibnamefont {Rasmussen}}, \ and\ \bibinfo {author} {\bibfnamefont {A.}~\bibnamefont {Veneziani}},\ }\bibfield  {title} {\enquote {\bibinfo {title} {Towards a theory for diffusive coupling functions allowing persistent synchronization},}\ }\href {\doibase 10.1088/0951-7715/27/3/501} {\bibfield  {journal} {\bibinfo  {journal} {Nonlinearity}\ }\textbf {\bibinfo {volume} {27}},\ \bibinfo {pages} {501} (\bibinfo {year} {2014})}\BibitemShut {NoStop}%
\bibitem [{\citenamefont {Cho}, \citenamefont {Nishikawa},\ and\ \citenamefont {Motter}(2017)}]{Cho_2017}%
  \BibitemOpen
  \bibfield  {author} {\bibinfo {author} {\bibfnamefont {Y.~S.}\ \bibnamefont {Cho}}, \bibinfo {author} {\bibfnamefont {T.}~\bibnamefont {Nishikawa}}, \ and\ \bibinfo {author} {\bibfnamefont {A.~E.}\ \bibnamefont {Motter}},\ }\bibfield  {title} {\enquote {\bibinfo {title} {Stable chimeras and independently synchronizable clusters},}\ }\href {\doibase 10.1103/PhysRevLett.119.084101} {\bibfield  {journal} {\bibinfo  {journal} {Phys. Rev. Lett.}\ }\textbf {\bibinfo {volume} {119}},\ \bibinfo {pages} {084101} (\bibinfo {year} {2017})}\BibitemShut {NoStop}%
\bibitem [{\citenamefont {Zhang}\ and\ \citenamefont {Motter}(2020)}]{Zhang_2020}%
  \BibitemOpen
  \bibfield  {author} {\bibinfo {author} {\bibfnamefont {Y.}~\bibnamefont {Zhang}}\ and\ \bibinfo {author} {\bibfnamefont {A.~E.}\ \bibnamefont {Motter}},\ }\bibfield  {title} {\enquote {\bibinfo {title} {Symmetry-independent stability analysis of synchronization patterns},}\ }\href {\doibase 10.1137/19M127358X} {\bibfield  {journal} {\bibinfo  {journal} {SIAM Review}\ }\textbf {\bibinfo {volume} {62}},\ \bibinfo {pages} {817--836} (\bibinfo {year} {2020})},\ \Eprint {http://arxiv.org/abs/https://doi.org/10.1137/19M127358X} {https://doi.org/10.1137/19M127358X} \BibitemShut {NoStop}%
\bibitem [{\citenamefont {Gambuzza}\ and\ \citenamefont {Frasca}(2019)}]{Gambuzza_2019}%
  \BibitemOpen
  \bibfield  {author} {\bibinfo {author} {\bibfnamefont {L.~V.}\ \bibnamefont {Gambuzza}}\ and\ \bibinfo {author} {\bibfnamefont {M.}~\bibnamefont {Frasca}},\ }\bibfield  {title} {\enquote {\bibinfo {title} {A criterion for stability of cluster synchronization in networks with external equitable partitions},}\ }\href {\doibase https://doi.org/10.1016/j.automatica.2018.11.026} {\bibfield  {journal} {\bibinfo  {journal} {Automatica}\ }\textbf {\bibinfo {volume} {100}},\ \bibinfo {pages} {212--218} (\bibinfo {year} {2019})}\BibitemShut {NoStop}%
\bibitem [{\citenamefont {Brady}, \citenamefont {Zhang},\ and\ \citenamefont {Motter}(2021)}]{Brady2021forget}%
  \BibitemOpen
  \bibfield  {author} {\bibinfo {author} {\bibfnamefont {F.~M.}\ \bibnamefont {Brady}}, \bibinfo {author} {\bibfnamefont {Y.}~\bibnamefont {Zhang}}, \ and\ \bibinfo {author} {\bibfnamefont {A.~E.}\ \bibnamefont {Motter}},\ }\href@noop {} {\enquote {\bibinfo {title} {Forget partitions: Cluster synchronization in directed networks generate hierarchies},}\ } (\bibinfo {year} {2021}),\ \Eprint {http://arxiv.org/abs/2106.13220} {arXiv:2106.13220 [nlin.AO]} \BibitemShut {NoStop}%
\bibitem [{\citenamefont {Eckmann}\ and\ \citenamefont {Ruelle}(1985)}]{Eckmann_1985}%
  \BibitemOpen
  \bibfield  {author} {\bibinfo {author} {\bibfnamefont {J.~P.}\ \bibnamefont {Eckmann}}\ and\ \bibinfo {author} {\bibfnamefont {D.}~\bibnamefont {Ruelle}},\ }\bibfield  {title} {\enquote {\bibinfo {title} {Ergodic theory of chaos and strange attractors},}\ }\href {\doibase 10.1103/RevModPhys.57.617} {\bibfield  {journal} {\bibinfo  {journal} {Rev. Mod. Phys.}\ }\textbf {\bibinfo {volume} {57}},\ \bibinfo {pages} {617--656} (\bibinfo {year} {1985})}\BibitemShut {NoStop}%
\bibitem [{\citenamefont {Pikovsky}\ and\ \citenamefont {Politi}(2016)}]{Pikovsky_politi_2016}%
  \BibitemOpen
  \bibfield  {author} {\bibinfo {author} {\bibfnamefont {A.}~\bibnamefont {Pikovsky}}\ and\ \bibinfo {author} {\bibfnamefont {A.}~\bibnamefont {Politi}},\ }\href {\doibase 10.1017/CBO9781139343473} {\emph {\bibinfo {title} {Lyapunov Exponents: A Tool to Explore Complex Dynamics}}}\ (\bibinfo  {publisher} {Cambridge University Press},\ \bibinfo {year} {2016})\BibitemShut {NoStop}%
\bibitem [{\citenamefont {Hindmarsh}\ and\ \citenamefont {Rose}(1984)}]{Hidmarsh_Rose1984}%
  \BibitemOpen
  \bibfield  {author} {\bibinfo {author} {\bibfnamefont {J.~L.}\ \bibnamefont {Hindmarsh}}\ and\ \bibinfo {author} {\bibfnamefont {R.~M.}\ \bibnamefont {Rose}},\ }\bibfield  {title} {\enquote {\bibinfo {title} {A model of neuronal bursting using three coupled first order differential equations},}\ }\href@noop {} {\bibfield  {journal} {\bibinfo  {journal} {Proc. R. Soc. Lond.}\ }\textbf {\bibinfo {volume} {221}},\ \bibinfo {pages} {87--102} (\bibinfo {year} {1984})}\BibitemShut {NoStop}%
\bibitem [{\citenamefont {Huang}\ \emph {et~al.}(2009)\citenamefont {Huang}, \citenamefont {Chen}, \citenamefont {Lai},\ and\ \citenamefont {Pecora}}]{huang2009generic}%
  \BibitemOpen
  \bibfield  {author} {\bibinfo {author} {\bibfnamefont {L.}~\bibnamefont {Huang}}, \bibinfo {author} {\bibfnamefont {Q.}~\bibnamefont {Chen}}, \bibinfo {author} {\bibfnamefont {Y.-C.}\ \bibnamefont {Lai}}, \ and\ \bibinfo {author} {\bibfnamefont {L.~M.}\ \bibnamefont {Pecora}},\ }\bibfield  {title} {\enquote {\bibinfo {title} {Generic behavior of master-stability functions in coupled nonlinear dynamical systems},}\ }\href@noop {} {\bibfield  {journal} {\bibinfo  {journal} {Physical Review E}\ }\textbf {\bibinfo {volume} {80}},\ \bibinfo {pages} {036204} (\bibinfo {year} {2009})}\BibitemShut {NoStop}%
\bibitem [{\citenamefont {Majhi}, \citenamefont {Perc},\ and\ \citenamefont {Ghosh}(2016)}]{Majhi2016}%
  \BibitemOpen
  \bibfield  {author} {\bibinfo {author} {\bibfnamefont {S.}~\bibnamefont {Majhi}}, \bibinfo {author} {\bibfnamefont {M.}~\bibnamefont {Perc}}, \ and\ \bibinfo {author} {\bibfnamefont {D.}~\bibnamefont {Ghosh}},\ }\bibfield  {title} {\enquote {\bibinfo {title} {Chimera states in uncoupled neurons induced by a multilayer structure},}\ }\href {\doibase 10.1038/srep39033} {\bibfield  {journal} {\bibinfo  {journal} {Scientific Reports}\ }\textbf {\bibinfo {volume} {6}},\ \bibinfo {pages} {39033} (\bibinfo {year} {2016})}\BibitemShut {NoStop}%
\bibitem [{\citenamefont {Majhi}, \citenamefont {Perc},\ and\ \citenamefont {Ghosh}(2017)}]{Majhi_2017}%
  \BibitemOpen
  \bibfield  {author} {\bibinfo {author} {\bibfnamefont {S.}~\bibnamefont {Majhi}}, \bibinfo {author} {\bibfnamefont {M.}~\bibnamefont {Perc}}, \ and\ \bibinfo {author} {\bibfnamefont {D.}~\bibnamefont {Ghosh}},\ }\bibfield  {title} {\enquote {\bibinfo {title} {{Chimera states in a multilayer network of coupled and uncoupled neurons}},}\ }\href {\doibase 10.1063/1.4993836} {\bibfield  {journal} {\bibinfo  {journal} {Chaos: An Interdisciplinary Journal of Nonlinear Science}\ }\textbf {\bibinfo {volume} {27}},\ \bibinfo {pages} {073109} (\bibinfo {year} {2017})},\ \Eprint {http://arxiv.org/abs/https://pubs.aip.org/aip/cha/article-pdf/doi/10.1063/1.4993836/14612458/073109\_1\_online.pdf} {https://pubs.aip.org/aip/cha/article-pdf/doi/10.1063/1.4993836/14612458/073109\_1\_online.pdf} \BibitemShut {NoStop}%
\bibitem [{\citenamefont {Omelchenko}\ \emph {et~al.}(2019)\citenamefont {Omelchenko}, \citenamefont {H{\"u}lser}, \citenamefont {Zakharova},\ and\ \citenamefont {Sch{\"o}ll}}]{omelchenko2019control}%
  \BibitemOpen
  \bibfield  {author} {\bibinfo {author} {\bibfnamefont {I.}~\bibnamefont {Omelchenko}}, \bibinfo {author} {\bibfnamefont {T.}~\bibnamefont {H{\"u}lser}}, \bibinfo {author} {\bibfnamefont {A.}~\bibnamefont {Zakharova}}, \ and\ \bibinfo {author} {\bibfnamefont {E.}~\bibnamefont {Sch{\"o}ll}},\ }\bibfield  {title} {\enquote {\bibinfo {title} {Control of chimera states in multilayer networks},}\ }\href@noop {} {\bibfield  {journal} {\bibinfo  {journal} {Frontiers in Applied Mathematics and Statistics}\ }\textbf {\bibinfo {volume} {4}},\ \bibinfo {pages} {67} (\bibinfo {year} {2019})}\BibitemShut {NoStop}%
\bibitem [{\citenamefont {Mikhaylenko}\ \emph {et~al.}(2019)\citenamefont {Mikhaylenko}, \citenamefont {Ramlow}, \citenamefont {Jalan},\ and\ \citenamefont {Zakharova}}]{mikhaylenko2019weak}%
  \BibitemOpen
  \bibfield  {author} {\bibinfo {author} {\bibfnamefont {M.}~\bibnamefont {Mikhaylenko}}, \bibinfo {author} {\bibfnamefont {L.}~\bibnamefont {Ramlow}}, \bibinfo {author} {\bibfnamefont {S.}~\bibnamefont {Jalan}}, \ and\ \bibinfo {author} {\bibfnamefont {A.}~\bibnamefont {Zakharova}},\ }\bibfield  {title} {\enquote {\bibinfo {title} {Weak multiplexing in neural networks: Switching between chimera and solitary states},}\ }\href@noop {} {\bibfield  {journal} {\bibinfo  {journal} {Chaos: an interdisciplinary journal of nonlinear science}\ }\textbf {\bibinfo {volume} {29}} (\bibinfo {year} {2019})}\BibitemShut {NoStop}%
\bibitem [{\citenamefont {Zakharova}(2020)}]{zakharova2020chimera}%
  \BibitemOpen
  \bibfield  {author} {\bibinfo {author} {\bibfnamefont {A.}~\bibnamefont {Zakharova}},\ }\bibfield  {title} {\enquote {\bibinfo {title} {Chimera patterns in networks},}\ }\href@noop {} {\bibfield  {journal} {\bibinfo  {journal} {Springer}\ } (\bibinfo {year} {2020})}\BibitemShut {NoStop}%
\bibitem [{\citenamefont {Blaha}\ \emph {et~al.}(2019)\citenamefont {Blaha}, \citenamefont {Huang}, \citenamefont {Della~Rossa}, \citenamefont {Pecora}, \citenamefont {Hossein-Zadeh},\ and\ \citenamefont {Sorrentino}}]{blaha2019cluster}%
  \BibitemOpen
  \bibfield  {author} {\bibinfo {author} {\bibfnamefont {K.~A.}\ \bibnamefont {Blaha}}, \bibinfo {author} {\bibfnamefont {K.}~\bibnamefont {Huang}}, \bibinfo {author} {\bibfnamefont {F.}~\bibnamefont {Della~Rossa}}, \bibinfo {author} {\bibfnamefont {L.}~\bibnamefont {Pecora}}, \bibinfo {author} {\bibfnamefont {M.}~\bibnamefont {Hossein-Zadeh}}, \ and\ \bibinfo {author} {\bibfnamefont {F.}~\bibnamefont {Sorrentino}},\ }\bibfield  {title} {\enquote {\bibinfo {title} {Cluster synchronization in multilayer networks: A fully analog experiment with l c oscillators with physically dissimilar coupling},}\ }\href@noop {} {\bibfield  {journal} {\bibinfo  {journal} {Physical review letters}\ }\textbf {\bibinfo {volume} {122}},\ \bibinfo {pages} {014101} (\bibinfo {year} {2019})}\BibitemShut {NoStop}%
\bibitem [{\citenamefont {Sorrentino}, \citenamefont {Pecora},\ and\ \citenamefont {Trajkovi{\'c}}(2020)}]{sorrentino2020group}%
  \BibitemOpen
  \bibfield  {author} {\bibinfo {author} {\bibfnamefont {F.}~\bibnamefont {Sorrentino}}, \bibinfo {author} {\bibfnamefont {L.}~\bibnamefont {Pecora}}, \ and\ \bibinfo {author} {\bibfnamefont {L.}~\bibnamefont {Trajkovi{\'c}}},\ }\bibfield  {title} {\enquote {\bibinfo {title} {Group consensus in multilayer networks},}\ }\href@noop {} {\bibfield  {journal} {\bibinfo  {journal} {IEEE Transactions on Network Science and Engineering}\ }\textbf {\bibinfo {volume} {7}},\ \bibinfo {pages} {2016--2026} (\bibinfo {year} {2020})}\BibitemShut {NoStop}%
\bibitem [{\citenamefont {Jiruska}\ \emph {et~al.}(2013)\citenamefont {Jiruska}, \citenamefont {de~Curtis}, \citenamefont {Jefferys}, \citenamefont {Schevon}, \citenamefont {Schiff},\ and\ \citenamefont {Schindler}}]{Jiruska2013epilepsy}%
  \BibitemOpen
  \bibfield  {author} {\bibinfo {author} {\bibfnamefont {P.}~\bibnamefont {Jiruska}}, \bibinfo {author} {\bibfnamefont {M.}~\bibnamefont {de~Curtis}}, \bibinfo {author} {\bibfnamefont {J.~G.~R.}\ \bibnamefont {Jefferys}}, \bibinfo {author} {\bibfnamefont {C.~A.}\ \bibnamefont {Schevon}}, \bibinfo {author} {\bibfnamefont {S.~J.}\ \bibnamefont {Schiff}}, \ and\ \bibinfo {author} {\bibfnamefont {K.}~\bibnamefont {Schindler}},\ }\bibfield  {title} {\enquote {\bibinfo {title} {Synchronization and desynchronization in epilepsy: controversies and hypotheses},}\ }\href {\doibase https://doi.org/10.1113/jphysiol.2012.239590} {\bibfield  {journal} {\bibinfo  {journal} {The Journal of Physiology}\ }\textbf {\bibinfo {volume} {591}},\ \bibinfo {pages} {787--797} (\bibinfo {year} {2013})},\ \Eprint {http://arxiv.org/abs/https://physoc.onlinelibrary.wiley.com/doi/pdf/10.1113/jphysiol.2012.239590} {https://physoc.onlinelibrary.wiley.com/doi/pdf/10.1113/jphysiol.2012.239590} \BibitemShut {NoStop}%
\bibitem [{\citenamefont {Pessoa}(2014)}]{pessoa2014understanding}%
  \BibitemOpen
  \bibfield  {author} {\bibinfo {author} {\bibfnamefont {L.}~\bibnamefont {Pessoa}},\ }\bibfield  {title} {\enquote {\bibinfo {title} {Understanding brain networks and brain organization},}\ }\href@noop {} {\bibfield  {journal} {\bibinfo  {journal} {Physics of life reviews}\ }\textbf {\bibinfo {volume} {11}},\ \bibinfo {pages} {400--435} (\bibinfo {year} {2014})}\BibitemShut {NoStop}%
\bibitem [{\citenamefont {Gambuzza}\ \emph {et~al.}(2021)\citenamefont {Gambuzza}, \citenamefont {Frasca}, \citenamefont {Sorrentino}, \citenamefont {Pecora},\ and\ \citenamefont {Boccaletti}}]{Gabuzza_2021}%
  \BibitemOpen
  \bibfield  {author} {\bibinfo {author} {\bibfnamefont {L.~V.}\ \bibnamefont {Gambuzza}}, \bibinfo {author} {\bibfnamefont {M.}~\bibnamefont {Frasca}}, \bibinfo {author} {\bibfnamefont {F.}~\bibnamefont {Sorrentino}}, \bibinfo {author} {\bibfnamefont {L.~M.}\ \bibnamefont {Pecora}}, \ and\ \bibinfo {author} {\bibfnamefont {S.}~\bibnamefont {Boccaletti}},\ }\bibfield  {title} {\enquote {\bibinfo {title} {Controlling symmetries and clustered dynamics of complex networks},}\ }\href {\doibase 10.1109/TNSE.2020.3037039} {\bibfield  {journal} {\bibinfo  {journal} {IEEE Transactions on Network Science and Engineering}\ }\textbf {\bibinfo {volume} {8}},\ \bibinfo {pages} {282--293} (\bibinfo {year} {2021})}\BibitemShut {NoStop}%
\bibitem [{\citenamefont {Bollt}\ \emph {et~al.}(2023)\citenamefont {Bollt}, \citenamefont {Fish}, \citenamefont {Kumar}, \citenamefont {Roque~dos Santos},\ and\ \citenamefont {Laurienti}}]{bollt2023fractal}%
  \BibitemOpen
  \bibfield  {author} {\bibinfo {author} {\bibfnamefont {E.}~\bibnamefont {Bollt}}, \bibinfo {author} {\bibfnamefont {J.}~\bibnamefont {Fish}}, \bibinfo {author} {\bibfnamefont {A.}~\bibnamefont {Kumar}}, \bibinfo {author} {\bibfnamefont {E.}~\bibnamefont {Roque~dos Santos}}, \ and\ \bibinfo {author} {\bibfnamefont {P.~J.}\ \bibnamefont {Laurienti}},\ }\bibfield  {title} {\enquote {\bibinfo {title} {Fractal basins as a mechanism for the nimble brain},}\ }\href@noop {} {\bibfield  {journal} {\bibinfo  {journal} {Scientific Reports}\ }\textbf {\bibinfo {volume} {13}},\ \bibinfo {pages} {20860} (\bibinfo {year} {2023})}\BibitemShut {NoStop}%
\bibitem [{\citenamefont {Siddique}\ \emph {et~al.}(2018)\citenamefont {Siddique}, \citenamefont {Pecora}, \citenamefont {Hart},\ and\ \citenamefont {Sorrentino}}]{Siddique2018}%
  \BibitemOpen
  \bibfield  {author} {\bibinfo {author} {\bibfnamefont {A.~B.}\ \bibnamefont {Siddique}}, \bibinfo {author} {\bibfnamefont {L.}~\bibnamefont {Pecora}}, \bibinfo {author} {\bibfnamefont {J.~D.}\ \bibnamefont {Hart}}, \ and\ \bibinfo {author} {\bibfnamefont {F.}~\bibnamefont {Sorrentino}},\ }\bibfield  {title} {\enquote {\bibinfo {title} {Symmetry- and input-cluster synchronization in networks},}\ }\href {\doibase 10.1103/PhysRevE.97.042217} {\bibfield  {journal} {\bibinfo  {journal} {Phys. Rev. E}\ }\textbf {\bibinfo {volume} {97}},\ \bibinfo {pages} {042217} (\bibinfo {year} {2018})}\BibitemShut {NoStop}%
\bibitem [{\citenamefont {Lin}\ \emph {et~al.}(2016)\citenamefont {Lin}, \citenamefont {Fan}, \citenamefont {Wang}, \citenamefont {Ying},\ and\ \citenamefont {Wang}}]{Lin_2016}%
  \BibitemOpen
  \bibfield  {author} {\bibinfo {author} {\bibfnamefont {W.}~\bibnamefont {Lin}}, \bibinfo {author} {\bibfnamefont {H.}~\bibnamefont {Fan}}, \bibinfo {author} {\bibfnamefont {Y.}~\bibnamefont {Wang}}, \bibinfo {author} {\bibfnamefont {H.}~\bibnamefont {Ying}}, \ and\ \bibinfo {author} {\bibfnamefont {X.}~\bibnamefont {Wang}},\ }\bibfield  {title} {\enquote {\bibinfo {title} {Controlling synchronous patterns in complex networks},}\ }\href {\doibase 10.1103/PhysRevE.93.042209} {\bibfield  {journal} {\bibinfo  {journal} {Phys. Rev. E}\ }\textbf {\bibinfo {volume} {93}},\ \bibinfo {pages} {042209} (\bibinfo {year} {2016})}\BibitemShut {NoStop}%
\end{thebibliography}%

\end{document}